\documentclass{amsart}
\usepackage{amsmath,amsfonts,amssymb,amsthm}

\makeatletter
\def\LaTeX{\leavevmode L\raise.42ex
    \hbox{\kern-.3em\size{\sf@size}{0pt}\selectfont A}\kern-.15em\TeX}
\makeatother

\newcommand{\BibTeX}{{\rm B\kern-.05em{\sc
          i\kern-.025emb}\kern-.08em\TeX}}

\makeatletter
\def\@currentlabel{2.1}\label{e:dispaa}
\def\@currentlabel{2.21}\label{e:dispau}
\def\@currentlabel{2.22}\label{e:dispav}
\def\@currentlabel{2.23}\label{e:dispaw}
\def\@currentlabel{2.24}\label{e:dispax}
\def\theequation{\thesection.\@arabic\c@equation}
\makeatother

\usepackage{afterpage}
\usepackage[dvips]{graphicx}
\usepackage{latexsym}

\newcommand{\cal}{\mathcal}
\newcommand{\dis}{\displaystyle}
 \newcommand{\epi}{{\bf P}}
 \newcommand{\epii}{{ P_i}}
 
 \newcommand{\epik}{{ P_k}}
 \newcommand{\pii}{{P_i}}
 
\newcommand{\epij}{{ P_j}}
\newcommand{\fr}{\frac{1}{2}}

\newcommand{\R}{{\mathbb{R}}}

\newcommand{\intr}{\int_{\R^N}}
\newcommand{\e }{\varepsilon }
\renewcommand{\theequation}{\arabic{section}.\arabic{equation}}

\DeclareMathOperator{\diag}{diag}

 \DeclareMathOperator{\spann}{span}
\newtheorem{theorem}{Theorem}[section]
\newtheorem{proposition}[theorem]{Proposition}
\newtheorem{lemma}[theorem]{Lemma}

\newtheorem{remark}[theorem]{Remark}

\begin{document}

\title[Clusters around saddle points]
{Positive and sign-changing clusters around saddle points\\ of the potential
for nonlinear elliptic problems}

\author{Teresa D'Aprile \and David Ruiz}

\address{Dipartimento di Matematica, Universit\`a di Roma
``Tor Vergata", via della Ricerca Scientifica 1, 00133 Roma,
ITALY, and Departamento de An{\'a}lisis Matem{\'a}tico, Universidad de
Granada, 18071 Granada, SPAIN.}

\thanks{T.D. has been supported by  the Italian PRIN Research Project 2007 \textit{Metodi
variazionali e topologici nello studio di fenomeni non lineari}}
\thanks{D.R. has been supported by the Spanish
Ministry of Science and Innovation under Grant MTM2008-00988 and by J. Andaluc\'{\i}a (FQM 116).}

\email{daprile@mat.uniroma2.it, daruiz@ugr.es}

\maketitle

\begin{abstract}

We study the existence and asymptotic behavior of positive and
sign-changing multipeak solutions for the equation
$$-\e^2\Delta v+V(x)v=f(v)\;\hbox{ in }\R^N,$$
where $\e$ is a small positive parameter,  $f$ is a superlinear,
subcritical and odd nonlinearity, $V$ is a uniformly positive
potential. No symmetry on $V$ is assumed. It is known (\cite{kw})
that this equation has positive multipeak solutions  with all
peaks approaching a local maximum of $V$. It is also proved that
solutions alternating positive and negative spikes exist in the
case of a minimum (see \cite{dapi2}). The aim of this paper is to
show the existence of both positive and sign-changing multipeak
solutions around a nondegenerate saddle point of $V$.
\bigskip



\end{abstract}

 \baselineskip 16pt

\section{Introduction}
\setcounter{equation}{0} This paper deals with the following
nonlinear perturbed elliptic equation
\begin{equation}\label{eq1}-\e^2\Delta v+V(x)v=|v|^{p-2}v\;\;\hbox{ in }\R^N\end{equation}
where  $N\geq 2$, $\e$ is a small parameter, the potential $V\in
{\mathcal C}^1(\R^N,\R)$
 is bounded from below away from zero, the exponent $p$ satisfies  $2<p<\frac{2N}{N-2}$ if
 $N\geq 3$ and $p>2$ if $N=2$. This equation arises when one looks for \textit{standing waves} of the
 nonlinear Schr\"odinger equation
$$i\e\frac{\partial \psi}{\partial t}=-\e^2\Delta \psi+V(x)\psi-|\psi|^{p-2}\psi,$$
which appears in different problems in nonlinear optics, in plasma
physics, etc.

Equation \eqref{eq1} has attracted much attention: a large number
of works are concerned with the question of semiclassical limit,
that is, the behaviour of solutions when $\e$ tends to zero. This
has an important physical interest since letting  $\e$ go to zero
formally describes the transition from Quantum Mechanics to
Classical Mechanics. It has been shown that if $P_0$ is a
nondegenerate or, more generally, a \textit{topologically
nontrivial} critical point of $V$, there exists a family of
solutions $v_\e$ which develops a single spike near $P_0$ as
$\e\to 0$ (\cite{abc}, \cite{ams}, \cite{gen2}, \cite{Del3},
\cite{grossi},
   \cite{yyl0},
 \cite{wang},  see \cite{amma} for further references). Also, when $V$ has several critical points,
  multi-peaks have been constructed with each peak concentrating  at a separate critical point
  (see \cite{Del2}, \cite{dfm}, \cite{gui}, \cite{oh} and references therein).

In this paper we are  interested in a special kind of solution for
equation \eqref{eq1}, the so called \textit{cluster}, i.e. a
combination of several interacting peaks concentrating at the same
point
 as $\e\to
0^+$. In \cite{kw} Kang and  Wei construct this kind of  solution:
more precisely, given $\ell\geq 1$ and $P_0$ a strict local
maximum of $V$, there exists a positive {\textit{cluster}} with
$\ell $ peaks concentrating at $P_0$. They also prove that such
solutions do not exist around nondegenerate minimum points of $V$.
After that, several papers have addressed the question of
existence of multibump solutions concentrating around a minimum of
$V$. This result has become known first for exactly one positive
and one negative peak (\cite{also}, \cite{bcw}), and later under
polygonal symmetries of $V$ (\cite{dapi1}), or in the
one-dimensional case (\cite{dft}). In a recent paper \textit{clusters}
with at most $6$ mixed positive and negative peaks have been
found, see \cite{dapi2}.

All previous results are concerned with the existence of a
\textit{clustered} solution localized around a minimum or a maximum point of
$V$. So the question of whether other critical points of $V$ may
generate a \textit{cluster} or not arises naturally. The aim of
this paper is to construct both positive and sign-changing
\textit{clusters} around a nondegenerate saddle point of $V$.

In order to provide the exact formulation of our results let us
fix some notation. We point out that most of the results contained
in the aforementioned papers can be extended to equations where
$|v|^{p-2}v$ is replaced by a more general nonlinear term. Then we
will consider the more general equation
\begin{equation}
\label{sch} \e^2\Delta v-V(x) v+f(v)=0\;\hbox{ in }\R^N
\end{equation} where we assume the following hypotheses on $f$:

\begin{enumerate}

\item[(f1)] $ f \in {\mathcal C}^{1+\sigma}_{loc} (\R)\cap
{\mathcal C}^2(0,+\infty)$ with $\sqrt{2}-1<\sigma<1 ;$
$f(0)=f'(0)=0$; $ f(t) =-f(-t)$ for all $t\in\R $; \item[(f2)]
$f(t)=O(t^{p-1})$ as $t\to +\infty$ for some $2<p<\frac{2N}{N-2}$ if
$N\geq 3$ and $p>2 $ if $N=2$;

\item[(f3)] the following \textit{limiting}  problem
\begin{equation}
\label{ground} \left\{
\begin{aligned}
 &\Delta w - w + f(w)=0, \;\;w>0 \mbox{ in } \R^N ,\\
  &w(0)=\max_{x \in \R^N} w(x),\;\; \lim_{ |x| \to +\infty} w(x)
  =0
  \end{aligned}
\right.
\end{equation} has a unique solution $w$, which is nondegenerate, i.e., denoting
by $L$ the linearized operator
$$L:H^2(\R^N)\to L^2(\R^N),\;\; L[u]:=\Delta u -u  + f' (w)u,$$
then
\begin{equation}
\label{Ker} \mbox{Kernel} (L) = \mbox{span} \left\{ \frac{\partial
w}{\partial x_1},\ldots, \frac{\partial w}{\partial x_N} \right\}.
\end{equation}
\end{enumerate}
By the well-known result of Gidas, Ni and Nirenberg (\cite{gnn})
$w$ is radially symmetric and strictly decreasing in $r=|x|$.
Moreover, by classical regularity results, the following
asymptotic behavior holds:
\begin{equation}
\label{wdecay}\begin{aligned} w(r), \;w''
(r)&=\frac{A}{r^{(N-1)/2}} e^{-r} \Big( 1+
O\Big(\frac{1}{r}\Big)\Big), \\ w'(r)&= -\frac{A}{r^{(N-1)/2}}
e^{-r} \Big( 1+ O\Big(\frac{1}{r}\Big)\Big),\end{aligned}
\end{equation}
where $ A>0$ is a suitable positive constant.

The class of nonlinearities $f$ satisfying (f1)-(f3) includes, and
it is not restricted to, the model $f(v)=|v|^{p-2} v$ with
$p>1+\sqrt{2}$ if $N=1,2$ and $p \in \left (1+ \sqrt{2},
\frac{2N}{N-2}\right )$ if $N\in [3,11]$ (if $N \geq 12$ the
interval is empty). Other nonlinearities can be found in
\cite{dancer}.

Let us now state the hypotheses on the potential $V$ that will be
used.

 \begin{enumerate}

\item[(V1)] $V\in {\mathcal C}^{2}(\R^N,\R)$ and $\inf_{\R^N}V>0$.

\item[(V2)] $V$ has a nondegenerate saddle point at $P_0$, and,
without loss of generality, we may assume $V(P_0)=1$. We define $r
\in \{ 1, \dots ,N-1\} $ as the number of positive eigenvalues of
$D^2V(P_0)$, counted with their multiplicity.

\end{enumerate}

As already mentioned, in this paper we give two results. First,
for any fixed positive integer $\ell$  there exists a $\ell$-peak
positive \textit{clustered} solution concentrating at $P_0$.
Furthermore each peak has a profile similar to $w$ suitably
rescaled. More precisely we will prove the following theorem.

\begin{theorem}
\label{th1} Assume that hypotheses (f1)--(f3) and (V1)--(V2) hold
and let $\ell\geq 1$ be fixed. Then, for $\e>0$ sufficiently
small, the equation (\ref{sch}) has a positive solution $v_\e\in
H^1(\R^N).$

\noindent Furthermore there exist $P_{1}^\e,\ldots,
P_\ell^\e\in\R^N$ such that, as $\e\to  0^+$,

\begin{enumerate}

\item[\rm(i)] $ v_\e (x) = \sum_{i=1}^\ell
w\big(\frac{x-P_{i}^\e}{\e} \big)+ o(\e)$ uniformly for
$x\in\R^N$;
\medskip

\item[\rm(ii)]  $|P_i^\e-P_j^\e| \geq 2 \beta \e\log \frac{1}{\e}$
($i\neq j$) and $|P_i^{\e}-P_0| \leq \e^{\beta}$ for any fixed
$\beta \in (0,1)$.

\end{enumerate}
\end{theorem}

Secondly, we prove that the equation  \eqref{sch} possesses a
{\textit{cluster}} with $h$ positive peaks and $k$ negative peaks
approaching $P_0$, where $h$ and $k$ are integers under some
restrictions. The exact formulation of the result is the
following. \noindent
\begin{theorem}
\label{th2} Assume that  $N \geq 2$ and hypotheses (f1)--(f3) and
(V1)--(V2)   hold.  Let $h,k$ satisfying $$h,\,k\geq 1,\quad
\ell:=h+k\leq 6.$$
\begin{enumerate}
\item[\rm(i)] If $r\geq 2$, then, for $\e>0$ sufficiently small,
the equation (\ref{sch}) has a solution $v_\e\in  H^1(\R^N).$

\noindent Furthermore there exist $P_{1}^\e,\ldots,
P_\ell^\e\in\R^N$ such that, as $\e\to  0^+$,

\begin{itemize}

\item $ v_\e (x) = \sum_{i=1}^h w\big(\frac{x-P_{i}^\e}{\e}
\big)-\sum_{i=h+1}^\ell w\big(\frac{x-P_{i}^\e}{\e} \big)+ o(\e)$
uniformly for $x\in\R^N$;
\medskip

\item  $|P_i^\e-P_j^\e| \geq 2 \beta \e\log \frac{1}{\e}$ ($i\neq
j$) and $|P_i^{\e}-P_0| \leq \e^{\beta}$ for any fixed $\beta \in
(0,1)$.

\end{itemize}
\item[\rm(ii)] If $r=1$, the same result as $\rm (i)$ holds with
the additional assumption $k\in\{h-1,h,h+1\}$.
\end{enumerate}

\end{theorem}

We point out that positive \textit{clustered} solutions have also been
found for the following equation
$$-\e^2\Delta v+v=Q(x)|v|^{p-2}v\quad \hbox{ in }\R^N$$
around critical points of  $Q$ (\cite{danyan5}, \cite{noya}). In
particular, as far as we know, the only work regarding
\textit{clusters} concentrating  near a saddle point  is
\cite{danyan5}. However we are unaware of \textit{cluster}
phenomena  with mixed positive and negative peaks near a saddle
point. Theorem \ref{th2} seems to be the first result in this
line.

The proofs of Theorem \ref{th1} and \ref{th2} rely on perturbation
arguments, which combine the variational approach with a
Lyapunov-Schmidt type procedure. A sketch of this procedure is
given in Section 2. Throughout the paper we will need some
asymptotic estimates, made in detail in Appendix A. With this
estimates in hand and thanks to the non-degeneracy condition
\eqref{Ker}, we can use the contraction mapping principle to solve
the auxiliary equation. Since the computations are quite
technical, they have been postponed to Appendix B.

In Section 3 we are concerned with the finite dimensional
bifurcation equation. Alternatively, we look for critical points
of an associated reduced functional. This is the main difficulty
of our problem; here the reduced functional has a quite involved
behaviour due to the different interactions of the potential and
the bumps.

It seems not easy to find the exact position of the bumps in a
direct way, if no symmetry assumptions are made. In this paper we
use a max-min technique applied to the reduced functional in the
spirit of \cite{dapi2}. This max-min argument is far from obvious,
specially in the case of sign-changing solutions. It takes into
account the interaction among the bumps (which depends on their
respective sign) and the effect of $V$ on each bump (which depends
on its spatial displacement).

\medskip

\noindent {\bf{NOTATION}}: Throughout the paper we will often use
the notation $C$ to denote generic positive constants.  The value
of $C$ is allowed to vary from place to place.

\section{The reduction process: sketch of the proof}
In this section we outline the main steps of the so called
{\textit{finite dimensional reduction}, which reduces the problem
to finding a critical point for a functional on a finite
dimensional space. We postpone the proofs and details to Appendix
A and Appendix B.

Associated to \eqref{sch} is the following energy functional:
   \begin{equation}\label{enfu}I_\e:H^1_{V}(\R^N)\to\R,\quad
   I_\e[v]:=\frac{1}{2} \int_{\R^N} \big(\e^2|\nabla v|^2+V(x)|v|^2\big) dx- \intr
   F(v) dx.\end{equation}
where $F(t)=\int_0^t f(s) ds$ and $$H_V^1(\R^N)=\Big\{v\in
H^1(\R^N)\,\Big|\, \intr V(x)|v|^2 dx <\infty\Big\}.$$
 Let us equip $H^1_V(\R^N)$ with the following scalar
product:
\begin{equation*}
\label{sca1} (u, v)_\e= \intr \big(\e^2\nabla u\nabla
v+V(x)uv\big)dx.
\end{equation*}
It is well known that $I_\e\in {\mathcal C}^2(H_V^1(\R^N),\R)$ and
the critical points of $I_\e$ are the finite-energy solutions of
\eqref{sch}.

Without loss of generality we assume throughout the paper that $P_0=0$. Moreover, after suitably rotating the
coordinate system, we may assume that in a small neighborhood  of $0$ the following expansion holds:
$$V(x)=1+\frac12\sum_{n=1}^N \lambda_n x_n^2+
o(|x|^2)\hbox{ as }x\to 0,$$ where $\lambda_n>0$ for
$n=1,\ldots,r$, $\lambda_{n}<0$ for $n=r+1,\ldots, N$.

Consider $M^+,\,M^-\in \R^{N^2}$ the following diagonal matrices
$$M^+=\diag (\lambda_1, \ldots, \lambda_r,0,\ldots, 0)\quad M^-=\diag (0,\ldots, 0,|\lambda_{r+1}|, \ldots, |\lambda_N|),$$
and set
$$M=M^+-M^-=D^2V(0)=\diag(\lambda_1,\ldots, \lambda_N),\quad \overline{M}=M^++ M^-=\diag (|\lambda_1|,\ldots,|\lambda_N|).$$
Next for  $\ell\geq 2$ define the configuration space:
\begin{equation*}\Gamma_{\e}=
\Bigg\{ {\bf P}=(P_1,\ldots,P_{\ell})\in\R^{N\ell}\,
\bigg|\,\overline{M}[P_i]^2< \e^{2\beta}\;\forall i,\;\;\;
w\Big(\frac{P_i-P_j}{\e}\Big)<\e^{2\beta}
 \hbox{ for }i\neq j\Bigg\},
\label{GA}
\end{equation*} where $\beta\in (\sigma,1)$ is a number sufficiently close to 1\footnote{
Observe that $\Gamma_\e$ is nonempty, since  for $\e$ sufficiently
small $\{{\bf P}\,|\,|P_i|\leq \e\log^2\frac{1}{\e},\,
|P_i-P_j|\geq 2\beta\e\log\frac{1}{\e}\hbox{ for }i\neq j\}\subset
\Gamma_\e$  thanks to  assumption (V2) and \eqref{wdecay}.}.
Observe that, according to \eqref{wdecay},
\begin{equation} \Gamma_\e\subset\Bigg\{ {\bf
P}=(P_1,\ldots,P_{\ell})\in\R^{N\ell}\, \bigg|\,|P_i|\leq
(\min_{i} |\lambda_i| )^{-1/2}\, \e^{\beta}\;\forall
i,\;\;\;|P_i-P_j|\geq 2\beta^2\e\log\frac{1}{\e}
 \hbox{ for }i\neq j\Bigg\}.
\label{GMC}
\end{equation}
For ${\bf P}=(P_1,\ldots,P_\ell)\in \Gamma_\e$ set
$$w_{P_i}(x)=w\Big(\frac{x-P_i}{\e}\Big), \quad
w_{\bf P}=\sum_{i=1}^\ell \tau_{i}
w_\epii,\;\;\;\tau_{i}\in\{-1,+1\}.$$ Let  $\chi\in \mathcal
C^\infty_0(\R^N)$ be a cut-off function such that $\chi(x)=1$ if
$|x|<1$, so that one has $\chi w_{ P_i},\,\chi w_\epi\in
H^1_V(\R^N)$.

We look for a solution to \eqref{sch} in a small neighbourhood of
the first approximation $\chi w_\epi$, i.e. a solution of the form
as $v:=\chi w_\epi+\phi,$ where the rest term $\phi$ is
\textit{small}. To this  aim we introduce the following functions:
\begin{equation*}
 Z_{P_i, n}= (V(x)-\e^2\Delta)\frac{\partial (\chi w_{\pii})}{\partial x_n}
 ,\quad i\in\{1,\ldots,\ell\},\; n\in\{1,\ldots,N\}.
\end{equation*}
The object is  to solve the following nonlinear problem: given
${\bf P }=(P_1, \ldots, P_\ell)\in \Gamma_\e$, find $(\phi,
\alpha_{in})$ such that
   \begin{equation}\label{nonl}\left\{\begin{aligned} & {\cal S}_{\e}[\chi w_\epi+\phi]=\sum_{i,n} \alpha_{in} Z_{P_i,n},\\
   &\phi\in H^2(\R^N)\cap H^1_V(\R^N),  \;\;\intr \phi Z_{P_i,n}\,dx  =0, \;\;\;i=1,\ldots, \ell
,\, n=1,\ldots,N,\end{aligned}\right.\end{equation} where
   \begin{equation}\label{operator}{\cal S}_{\e}[v]=\e^2\Delta v-V(x)v + f(v).\end{equation}

   \begin{lemma}\label{reg} Set $ \eta =\beta^2(1+\sigma)
$. Provided that $\e>0$ is sufficiently small,
   for every    $ {\bf P }\in \Gamma_\e$ there is a pair $ (\phi_{\bf P }, \alpha_{in}({\bf P})) \in
\big(H^2(\R^N)\cap H^1_V(\R^N)\big)\times\R^{N\ell} $ satisfying
(\ref{nonl}) and
\begin{equation}
   \label{phi}
 \|\phi_{\bf P }\|_\infty\leq C\e^\eta,\;\; (\phi_{\bf P }, \phi_{\bf P})_\e\leq C \e ^{N+2\eta} ,\;\;|\alpha_{in}({\bf P})|\leq C\e^{1+\eta}.
   \end{equation}
 Moreover
 the map ${\bf P}\in\Gamma_\e\mapsto \phi_{\bf P}\in H^1_V(\R^N)$ is ${\mathcal C}^1$.

   \end{lemma}

We refer to Appendix B for the proof.

 For $\e>0$ sufficiently small consider the reduced
functional
\begin{equation*}
\label{Mept}  J_\e: \Gamma_\e\to\R,\;\; J_\e[{\bf P } ]:=
\e^{-N}I_\e [ \chi w_{\epi} + \phi_{{\bf P }}] - c_1 ,
\end{equation*}
 where $ \phi_{\bf P }$ has been constructed in Lemma \ref{reg} and
 $c_1=\frac{\ell}{2}\intr \big(|\nabla  w|^2+w^2\big)dx-\ell\intr F(w) dx $.
The next proposition  contains the key expansions of $J_\e$ and
$\nabla J_\e$ (see Appendix B for the proof).
\begin{proposition}\label{energy111}
The following expansions hold:
\begin{equation}\label{preespi1} \begin{aligned}J_{\e} [{\bf P }]=&
\frac{c_2}{2}\sum_{i=1}^\ell M[P_i]^2-\frac{1}{2} \e^{-N} \sum_{i\neq
j}\tau_{i}\tau_{j} \int_{\R^N} f(w_{P_i})w_{P_j}
dx+o(\e^{2\beta}),
\end{aligned}\end{equation}
\begin{equation}\label{preespi2}\frac{\partial J_\e}{\partial P_{i}} [{\bf P}]=c_2M[P_i]-\e^{-N}
\sum_{j,\,j\neq i} \tau_{i}\tau_{j} \frac{\partial}{\partial
P_i}\bigg[\intr
  f(w_\epii)
w_\epij\,dx\bigg] +o(\e^{\beta}), \quad i=1,\ldots, \ell
\end{equation}
uniformly for ${\bf P }\in \Gamma_\e,$ where
 $c_2=\frac12\intr w^2 dx$.

By Lemma \ref{gamma0} (see Appendix A), we have also the following
expansion:
\begin{equation}\label{espi1} \begin{aligned}J_{\e} [{\bf P }]=&\frac{c_2}{2}\sum_{i=1}^\ell M[P_i]^2-\frac{c_3}{2}\sum_{i\neq j}\tau_{i}\tau_{j} w\Big(\frac{P_i-P_j}{\e}\Big)+o(\e^{2\beta}),
\end{aligned}\end{equation}
uniformly for ${\bf P }\in \Gamma_\e,$ where
 $c_3=\intr f(w) e^{x_1}dx$. \label{energy3}
\end{proposition}

\noindent Finally the next lemma  concerns the relation between
the critical points of $J_\e $ and those of $I_\e $. It is quite
standard in singular perturbation theory; its proof can be found
in \cite{amma}, for instance.

\begin{lemma}\label{relation}
Let ${\bf P }_\e  \in \Gamma_\e$ be a critical point of $J_\e $.
Then, provided that $\e>0$ is sufficiently small, the
corresponding function $ v_\e=\chi w_{\epi_\e} + \phi_{ {\bf P
}_\e}$ is  a solution of (\ref{sch}). \end{lemma}

So, we conclude the proof by showing the existence of a critical
point of $J_{\e}$. This will be accomplished in next section.

\section{A max-min argument: proof of Theorem \ref{th1} and Theorem \ref{th2}}
In this section we apply a max-min argument to characterize a
topologically nontrivial  critical value of $J_\e$. More precisely
we will construct sets $\mathcal D_\e$, $K$, $K_0\subset
\R^{N\ell}$ satisfying  the following properties:
\begin{enumerate}
\item [(P1)] $\mathcal D_\e$ is an open set,  $K_0$ and $K$ are
compact sets,  $K$ is connected and
\begin{equation*}K_0\subset K\subset \mathcal D_\e\subset {\overline {\mathcal D}}_\e\subset \Gamma_{\e};\end{equation*}
\item[(P2)] if  we define the complete metric space ${\mathcal F}$
by
 $${\mathcal
F}=\{\eta:K\to {\mathcal D}_\e\,|\, \eta\hbox{
continuous},\;\eta({\bf P})={\bf P}\;\forall {\bf P}\in K_0\},$$
then
\begin{equation}\label{mima}{ J}^*_\e:=\sup_{\eta\in{\mathcal F}}\min_{{\bf P}\in
K}J_\e[\eta({\bf P})]<\min_{{\bf P}\in K_0} J_\e[{\bf
P}].\end{equation} \item[(P3)] For every ${\bf
P}\in\partial\mathcal D_\e$ such that $J_\e[{\bf P}]={J}^*_\e$, we
have that $\partial \mathcal D_\e$ is smooth at ${\bf P}$ and
there exists a vector $\tau_{\bf P}$ tangent to $\partial\mathcal
D_\e$ at ${\bf P}$ so that $J_\e'[{\bf P}](\tau_{\bf P})\neq 0$.
\end{enumerate}

Under these assumptions a critical point ${{\bf P}}_\e\in
{\mathcal D}_\e$ of $J_\e$ with $J_\e[{{\bf P}}_\e]={J}^*_\e$
exists, as a standard deformation argument involving the gradient
flow of $J_\e$ shows.

 We define

  $${\mathcal D}_\e=\bigg\{{\bf P}\in \R^{N\ell}\;
\bigg|\,c_2\sum_{i=1}^\ell \overline{M}[P_i]^2+c_3\sum_{i\neq
j}w\Big(\frac{P_i-P_j}{\e}\Big) <c_4\e^{2\beta} \bigg\}.$$ where
$c_4=\min\{c_2,c_3\}$. We immediately get $\overline{{\mathcal
D}}_\e\subset \Gamma_\e$. In the following we will denote by $A$
and $B$ the subspaces associated to the positive and negative
eigenvalues of $M$ respectively, whose direct sum is $\R^N$, i.e.
$$A=\spann\{{\bf e}_1,\ldots, {\bf e}_r\}, \quad B=\spann\{{\bf e}_{r+1},\ldots, {\bf e}_N\},$$ where ${\bf e}_1,\ldots, {\bf e}_N$ is the standard basis in $\R^N$.

\subsection{Definition of $K$, $K_0$, and proof of (P1)-(P2)}
In this subsection we define the sets $K,\,K_0$ for which
properties (P1)-(P2) hold. In addition, we will prove that
\begin{equation}\label{crucial}J_\e^*=o(\e^{2\beta}).\end{equation}

For the sake of clarity we distinguish  the case of positive peaks
from that of mixed positive and negative peaks.
\subsubsection{I case: $k=0$, $\ell=h\geq 1$}
We have $\tau_i=1$ for all  $i=1,\ldots,\ell$. Let us fix
$b_1,\ldots, b_\ell\in B$ such that
$$|b_i|\leq 2\ell \e \log\frac1\e,\,\forall i,\quad |b_i-b_j|\geq 2\e\log\frac1\e\hbox{ for }i\neq j,$$
and define the following convex
open  set $U$ of $A^\ell$:
 $$U=\bigg\{(a_1,\ldots, a_\ell)\in A^{\ell}\;\bigg|\;  c_2\sum_{i=1}^\ell M^+[a_i]^2<\frac{ c_4}{2}\e^{2\beta}\bigg\},$$
 and
$$K=\bigg\{{\bf P}=(a_1+b_1,\ldots, a_\ell+b_\ell)\in \R^{N\ell}\;\Big|\; (a_1,\ldots, a_\ell)\in\overline{U}\bigg\}, $$
$$K_0:=\bigg\{{\bf P}=(a_1+b_1,\ldots, a_\ell+b_\ell)\in \R^{N\ell}\;\bigg|\;
(a_1,\ldots, a_\ell)\in\partial U\}.$$
$K$ is clearly isomorphic to $\overline{U}$
by  the immediate isomorphism $$(a_1+b_1,\ldots, a_\ell+b_\ell)\in
K\longleftrightarrow (a_1,\ldots, a_\ell)\in\overline{U}$$
  and $K_0\approx \partial U$.
 $K_0$ and
$K$ are compact sets, $K$ is connected (since $\overline{U}$ is convex) and
$K_0 \subset K$. Furthermore
$M[a_i+b_i]^2=M^+[a_i]^2-M^-[b_i]^2=M^+[a_i]^2+O(\e^2\log^2\frac1\e)$, analogously $\overline{M}[a_i+b_i]^2=M^+[a_i]^2+M^-[b_i]^2=M^+[a_i]^2+O(\e^2\log^2\frac1\e)$
and, since $w$ is decreasing in $|x|$,
$w(\frac{a_i+b_i-a_j-b_j}{\e})\leq w(\frac{b_i-b_j}{\e})=o(\e^2)$
for $i\neq j$. Then we deduce $K\subset {\cal D}_\e$ and, by
Proposition \ref{energy111},
\begin{equation}\label{onkappa}J_\e({\bf P})=\frac{c_2}2\sum_{i=1}^\ell M^+[a_i]^2+o(\e^{2\beta})\hbox{ uniformly on }K,\end{equation}
 by which, since $c_2\sum_{i=1}^\ell M^+[a_i]^2=\frac{ c_4}{2}\e^{2\beta}$ if  $(a_1,\ldots, a_\ell)\in\partial U$,
\begin{equation}\label{onkappa0}J_\e[{\bf P}]=\frac{c_4}{4}(1+o(1))\e^{2\beta} \hbox{ uniformly on }K_0.\end{equation}
 Let $\eta\in {\mathcal F},$ namely
$\eta:K\to {\mathcal D}_\e$ is a continuous function such that
$\eta({\bf P})={\bf P}$ for any $ {\bf P}\in K_0.$
 Then we can compose the following maps

 $$A^\ell\supset \overline{U}\longleftrightarrow K\stackrel{\eta}{\longrightarrow}
 \eta(K)\subset {\cal D}_\e\stackrel{ (\pi_A)^\ell}{\longrightarrow} {A^\ell},$$
 denoting by $\pi_A$ the orthogonal projection of $\R^N$ onto $A$,
 and we call $T:\overline{U}\to A^\ell$ the resulting composition. $T$ is a continuous map. We claim that $T=id$ on $\partial U$. Indeed, if $(a_1,\ldots, a_\ell)\in \partial U$, then  $(a_1+b_1,\ldots, a_\ell+b_\ell)\in K_0$, consequently $\eta (a_1+b_1,\ldots, a_\ell+b_\ell)=a_1+b_1,\ldots, a_\ell+b_\ell$, by which
 $$T(a_1,\ldots, a_\ell)=(\pi_A)^\ell(a_1+b_1,\ldots, a_\ell+b_\ell)=(a_1,\ldots, a_\ell).$$
 Since $0=(0,\ldots, 0)\in U$,  hence the theory of the topological degree ensures
 that $\deg(T,U,0)=\deg(id, U,0)=1$. Then there exists $(\bar{a}_1,\ldots, \bar{a}_\ell)\in U$ such that $T(\bar{a}_1,\ldots, \bar{a}_\ell)=0$,
 i.e. ${\bf P}^\eta:=\eta(\bar{a}_1+b_1, \ldots, \bar{a}_\ell+b_\ell)\in B^\ell$. Using Proposition \ref{energy111} we get
 $$\min_{{\bf P}\in K}J_\e[\eta({\bf P})]\leq J_\e[{\bf P}^\eta]=
 - \frac{c_2}{2}\sum_{i=1}^\ell M^-[P^\eta_i]^2-\frac{c_3}{2}\sum_{i\neq j}w\Big(\frac{P_i^\eta-P_j^\eta}{\e}\Big)+o(\e^{2\beta})\leq o(\e^{2\beta}).$$
Hence
$$J_\e^*=\sup_{\eta\in{\cal F}}\min_{{\bf P}\in K}J_\e[\eta({\bf P})]\leq o(\e^{2\beta}).$$
On the other hand, by taking $\eta=id$ and using \eqref{onkappa},
 $$J_\e^*\geq \min_{{\bf P}\in K}J_\e[{\bf P}]\geq o(\e^{2\beta}).$$
Combining last two estimates we get \eqref{crucial}. Finally
comparing \eqref{crucial} with \eqref{onkappa0}, the max-min
inequality \eqref{mima} follows.

\subsubsection{II case: $h, \,k\geq 1,\, \ell:h+k\leq 6$}
For the sake of simplicity assume $h\geq k$.  In order to define
$K$ and $K_0$, we need to consider special type of configurations
${\bf P}(a, {\bf r})\in\R^{N\ell} $. To this aim it is convenient
to  distinguish three cases.

\begin{enumerate}
\item[(i)] \fbox{$k=h$ or $k=h-1$.}

We set $\tau_i=(-1)^{i+1}$. Let us fix ${\bf v}\in A$ such that
$|{\bf v}|=1$. Then we consider  the configurations  ${\bf P}$
which lie  in $A$  and  are aligned in the direction ${\bf v}$
with alternating sign, i.e. configurations of the form
\begin{equation}\label{firstform}{\bf P}={\bf P}(a,{\bf r})=\left(\begin{array}{c}a\\ a+r_2 {\bf v}\\  a+(r_2+r_3){\bf v}\\  \ldots\\  a+(r_2+\ldots
+r_\ell){\bf v}
\end{array}\right)\in \R^{N\ell},\end{equation}
where $a\in A$ and ${\bf r}=(r_2,\ldots, r_\ell)\in
(0,+\infty)^{\ell-1}$. Observe that by construction we have
$|P_i-P_j|=r_{j+1}+r_{j+2}+\ldots+r_i$ if $i>j$, therefore
\begin{equation}\label{firstform1}r_i=|P_i-P_{i-1}|=\min_{j<i}|P_i-P_j|
=\min_{j<i, \,\tau_j=-\tau_i}|P_i-P_j|\quad \forall i=2,\ldots,
\ell.\end{equation} Moreover, if $i>j$ and $\tau_i\tau_j=1$, then
$i\geq j+2$, and consequently $|P_i-P_j|\geq r_i+r_{i-1}$, by
which
\begin{equation}\label{firstform2}|P_i-P_j|\geq 2\min_{2\leq s\leq\ell}r_s\;\hbox{ if }\tau_i\tau_j=1, i\neq j.\end{equation}
\item[(ii)]  \fbox{$k=1,$ $2\leq h \leq 5$, $r\geq 2$}

We set $\tau_1=-1$, $\tau_i=1$ for $i=2,\ldots, \ell$. Let us fix
${\bf v}_2,\ldots,{\bf v}_\ell\in A$ such that $|{\bf v}_i|=1$ and
each ${\bf v}_i$  points at the vertex of a regular $h$ polygon.
Then we consider the configurations ${\bf P}$ which lie in  $A$
and such that $P_i$ for $i\geq 2$ is located on the half-line
starting from $P_1$   in the ${\bf v}_i$ direction; more precisely
\begin{equation}\label{secondform}{\bf P}={\bf P}(a,{\bf r})=\left(\begin{array}{c}a\\ a+r_2 {\bf v}_2\\  a+r_3{\bf v}_3\\
\ldots\\  a+r_\ell{\bf v}_\ell
\end{array}\right)\in \R^{N\ell}, \end{equation}where $a\in A$ and ${\bf r}=(r_2,\ldots, r_\ell)\in (0,+\infty)^{\ell-1}.$
We point out that
\begin{equation}\label{secondform1}r_i=|P_i-P_1|=\min_{j<i,\,
\tau_j=-\tau_i}|P_i-P_j|\quad \forall i=2,\ldots,
\ell.\end{equation} Moreover, if $\tau_i\tau_j=1$, then $i,j\geq
2$  and we have $|P_i-P_j|^2=
r_i^2+r_j^2-2r_ir_j\cos\frac{2\pi(i-j)}{h}$, which implies
 \begin{equation}\label{secondform2}|P_i-P_j|\geq\sqrt{2-2\cos\frac{2\pi}{h}}\min_{2\leq s\leq  \ell}r_s\;\hbox{ if }\tau_i\tau_j=1, \;i\neq j.
\end{equation}
Taking into account that $\sqrt{2-2\cos\frac{2\pi}{h}}>1$  if $2\leq h\leq 5$, by \eqref{secondform1} and \eqref{secondform2} we
immediately get
 \begin{equation}\label{secondform3}
 |P_i-P_j|\geq\min_{2\leq s\leq  \ell}r_s\hbox{ if }i\neq j.
 \end{equation}

 \item[(iii)] \fbox{$h=4$, $k=2$, $r\geq 2$}

  We set $\tau_1=\tau_3=-1$, $\tau_2=\tau_4=\tau_5=\tau_6=1$. Let us
 fix two orthogonal  vectors ${\bf v}, \, {\bf w}\in A$ such that $|{\bf v}|=|{\bf w}|=1$. Then consider  the configurations of the type
\begin{equation}\label{thirdform} {\bf P}(a,{\bf r})=\left(\begin{array}{c}a\\ a+r_2 {\bf v}\\  a+(r_2+r_3){\bf v}\\  a+(r_2+r_3+r_4){\bf v}\\
a+r_5{\bf w}\\ a-r_6 {\bf w}
\end{array}\right)\in \R^{N\ell}, \end{equation} where $a\in A$ and $ {\bf r}=(r_2,\ldots, r_\ell)\in (0,+\infty)^{\ell-1}.$
It is immediate to check that
\begin{equation}\label{thirdform1}r_i=\min_{j<i}|P_i-P_j|=\min_{j<i,\,\tau_j=-\tau_i}|P_i-P_j|\end{equation}
and
\begin{equation}\label{thirdform2}|P_i-P_j|\geq \sqrt{2}\min_{2\leq s\leq \ell}r_s\;\hbox{ if }\tau_i\tau_j=1, \;i\neq j.\end{equation}
\end{enumerate}
\bigskip

Observe that (i)-(ii)-(iii) cover all cases $(h,k)$ with the assumptions of Theorem \ref{th2}.

We now define:
$$ S=\left\{{\bf P} \in \R^{N\ell}\;\Bigg|\;
 c_2\sum_{i=1}^\ell
 \overline{M}[P_i]^2+c_3\sum_{i\neq j}w\Big(\frac{P_i-P_j}{\e}\Big) < \frac{c_4}{2}\e^{2\beta}\right\} $$
and
$$ \widetilde{U}= \left\{(a, {\bf r})\in A\times (0,+\infty)^{\ell-1}: \ {\bf P}(a, {\bf r}) \in S\right\} .$$
$\widetilde{U}$ is an open set. In principle, we do not know
whether $\widetilde{U}$ is connected or not, so we will define $U$
as a conveniently chosen connected component. We claim that
$(0,{\bf r}_\e) \in \widetilde{U} $, where ${\bf r}_{\e}$ is
defined as: ${\bf r}_\e=\big(2\e\log\frac 1\e,\ldots,
2\e\log\frac1\e\big)\in \R^{\ell-1}$.

Indeed, setting ${\bf P}(0,{\bf r}_\e)=(P_1^0,\ldots,
P_\ell^0)$, according to \eqref{firstform1}, \eqref{secondform3}
and \eqref{thirdform1} we have $|P_i^0-P_j^0|\geq 2\e\log\frac1\e$
for $i\neq j$ and, by \eqref{wdecay}, we immediately check
$w(\frac{P^0_i-P_j^0}{\e})=o(\e^2)$ for $i\neq j$. Furthermore,
according to
\eqref{firstform}-\eqref{secondform}-\eqref{thirdform} one has
$|P_i^0|\leq 2\ell\e\log\frac1\e$; then we infer $M[P^0_i]^2=
M^+[P^0_i]^2=O(\e^2\log^2\frac1\e)$.

Now we are in conditions of defining $U$, $K$ and $K_0$:
$$U=\mbox{ the connected component of } \widetilde{U} \mbox{ containing } (0,\bf{r}_\e), $$
$$ K=\left\{ {\bf P}(a, {\bf r}) \in \R^{N\ell}:\ (a, {\bf r}) \in \overline{U} \right\}, $$
$$ K_0=\left\{ {\bf P}(a, {\bf r}) \in \R^{N\ell}:\ (a, {\bf r}) \in \partial U \right\}.$$
$K$ is clearly isomorphic to $\overline{U}$ by the
obvious isomorphism, and $K_0 \approx \partial U$. In particular $K$ and $K_0$ are compact sets and $K$ is connected.
Moreover we have $K_0\subset K\subset {\cal D}_\e$.

If ${\bf P}={\bf P}(a,{\bf r})\in K$, by \eqref{GMC} and \eqref{firstform1}, \eqref{secondform1},
\eqref{thirdform1} we get $r_i\geq 2\beta^2\e\log\frac1\e$ for all $i=2, \ldots,\ell$; then \eqref{firstform2},
\eqref{secondform2}, \eqref{thirdform2} imply $|P_i-P_j|\geq 2\e\log\frac1\e$ if $\tau_i\tau_j=1$ and $i\neq j$
for $\beta$ sufficiently close to 1 (observe that $\sqrt{2-2\cos\frac{2\pi}{h}}>1$ if $2\leq h\leq 5$). Then  by
\eqref{wdecay} it follows that
$$w\Big(\frac{P_i-P_j}{\e}\Big)=o(\e^2) \;\hbox{ if }\tau_i\tau_j=1, \;i\neq j,\,\hbox{ uniformly on  }K.$$
Roughly speaking, the configurations in $K$
have the crucial property that  the mutual distance between
the points $P_i$, $P_j$ with $\tau_i=\tau_j$, $i\neq j$,  is
sufficiently large
so that their interaction term
$w(\frac{P_i-P_j}{\e })$ becomes negligible; moreover, since $K\subset A^\ell$, then $M[P_i]^2=M^+[P_i]^2$,   and
consequently  the main   terms which appear in  $J_\e$ are
positive. Indeed by Proposition \ref{energy111} we deduce
\begin{equation}\label{1onkappa}J_\e[{\bf
P}]=\frac{c_2}2\sum_{i=1}^\ell M^+[P_i]^2+\frac{c_3}{2}\sum_{i\neq
j}w\Big(\frac{P_i-P_j}{\e}\Big)+o(\e^{2\beta})\hbox{ uniformly on
}K,\end{equation}
 by which, since $c_2\sum_{i=1}^\ell M^+[P_i]^2+c_3\sum_{i\neq
j}w\big(\frac{P_i-P_j}{\e}\big)=\frac{c_4}{2}\e^{2\beta} $ if ${\bf P}\in K_0$,
\begin{equation}\label{1onkappa0}J_\e[{\bf P}]=\frac{c_4}{4}(1+o(1))\e^{2\beta} \hbox{ uniformly on }K_0.\end{equation}

Let $\eta\in {\mathcal F},$ namely $\eta:K\to {\mathcal D}_\e$ is
a continuous function such that $\eta({\bf P})={\bf P}$ for any $
{\bf P}\in K_0.$
 Then we can compose the following maps

 $$A\times (0,+\infty)^{\ell-1} \supset \overline{U}\longleftrightarrow K\stackrel{\eta}{\longrightarrow} \eta(K)\subset {\cal D}_\e\stackrel{ {\cal H}}{\longrightarrow} {A}\times (0,+\infty)^{\ell-1}$$  where ${\cal H}=({\cal H}_1, \ldots, {\cal H}_\ell):\R^{N\ell}\to A\times(0,+\infty)^{\ell-1}$ is defined by
 $${\cal H}_1(P_1,\ldots, P_\ell)=\pi_A(P_1), \quad {\cal H}_i(P_1,\ldots, P_\ell)=\min_{j<i,\,\tau_j=-\tau_i}|P_i-P_j|\,\hbox{ for }i\geq 2,$$
 denoting by $\pi_A$ the orthogonal projection of $\R^N$ onto $A$.
 We set $$T:\overline{U}\to A\times (0,+\infty)^{\ell-1}$$ the resulting composition.
 Clearly $T$ is a continuous map. We claim that $T=id$ on $\partial U$.
Indeed, if $(a,{\bf r})\in\partial U$, then by construction  ${\bf
P}(a,{\bf r})\in K_0$; consequently $\eta({\bf P}(a,{\bf r}))={\bf
P}(a,{\bf r})$, by which, using the definitions
\eqref{firstform}-\eqref{secondform}-\eqref{thirdform},
$${\cal H}_1({\bf P}(a,{\bf r}))=\pi_A(a)=a$$
while, using
\eqref{firstform1}-\eqref{secondform1}-\eqref{thirdform1},
$${\cal H}_i({\bf P}(a,{\bf r}))=r_i\;\hbox{ for }i\geq 2.$$
This proves that $T=id$ on $\partial U$.

The theory of the topological degree assures that $\deg(T, U,
(0,{\bf r}_\e))=\deg (id, U, (0,{\bf r}_\e))=1$; then there exists
$({a}_\eta,{\bf r}_\eta)\in U$ such that $T( a_\eta,{{\bf
r}}_\eta)= (0, {\bf r}_\e)$, i.e., setting   ${\bf
P}^\eta:=\eta({\bf P}({ a}_\eta,{{\bf r}}_\eta))\in\eta(K)$,
\begin{equation}\label{chat}\pi_A(P_1^\eta)=0,\quad \min_{j<i, \,\tau_j=-\tau_i}|P_j^\eta-P_i^\eta|=2\e\log\frac1\e\;\hbox{ for }i\geq 2.\end{equation}
In particular this implies $$P_1^\eta\in B,\quad \quad
|P_j^\eta-P_i^\eta|\geq 2\e\log\frac1\e \;\hbox{ if
}\tau_i=-\tau_j,$$
 which gives $$w\Big(\frac{P_i^\eta-P_j^\eta}{\e}\Big)=o(\e^2)\;\hbox{ if }\tau_i=-\tau_j.$$
Moreover, by the second of \eqref{chat}, recalling that
$\tau_1=-\tau_2$, it is not difficult to check that
$|P_i^\eta-P_1^\eta|\leq 2\ell\e\log\frac1\e $ for all $i$; then
$|\pi_A(P_i^\eta)|=|\pi_A(P_i^\eta-P_1^\eta)|\leq
2\ell\e\log\frac1\e$ and consequently
$$M[P_i^\eta]^2=-M^-[P_i^\eta]^2+O\Big(\e^2\log^2\frac1\e\Big).$$
By Proposition \ref{energy111} we infer $$\min_{{\bf P}\in
K}J_\e[\eta({\bf P})]\leq J_\e[{\bf P}^\eta]=
 - \frac{c_2}{2}\sum_{i=1}^\ell M^-[P^\eta_i]^2-\frac{c_3}{2}\sum_{i\neq j,\, \tau_i=\tau_j}w\Big(\frac{P_i^\eta-P_j^\eta}{\e}\Big)+o(\e^{2\beta})\leq o(\e^{2\beta}).$$
 By taking the supremum for all the maps $\eta\in{\cal F}$ we obtain  $$J_\e^*=\sup_{\eta\in{\cal F}}\min_{{\bf P}\in K}J_\e[\eta({\bf P})]\leq o(\e^{2\beta}).$$
On the other hand, by taking $\eta=id$ and using \eqref{1onkappa},
 $$J_\e^*\geq \min_{{\bf P}\in K}J_\e[{\bf P}]\geq o(\e^{2\beta}).$$
Last two estimates yield \eqref{crucial}. Finally, comparing
\eqref{crucial} with \eqref{1onkappa0}, the max-min inequality
\eqref{mima} follows.

\subsection{Proof of (P3)}

Let us define $$\Phi_{\e}: \Gamma_\e \to \R, \ \Phi_\e({\bf
P})= \frac{c_2}{2}\sum_{i=1}^\ell \overline{M}[P_i]^2+
\frac{c_3}{2}\sum_{i\neq j}w\Big(\frac{P_i-P_j}{\e}\Big).$$

We shall prove (P3) by contradiction: assume that there exist
$\e_n \to 0$, ${\bf P}_{\e_n}=(P^{\e_n}_1,\ldots,P^{\e_n}_\ell)
\in\partial\mathcal D_{\e_n}$ and a vector $(\mu_{\e_n, 1},\ \mu_{\e_n, 2})$ in the unit circle, i.e. $\mu_{\e_n, 1}^2+ \mu_{\e_n, 2}^2=1$, such that:
$$ \Phi_{\e_n} ({\bf P}_{\e_n}) = \frac{c_4}{2}\e_n^{2\beta}, $$
$$ J_{\e_n}[{\bf P}_{\e_n}]={J}^*_{\e_n}, $$
$$ \mu_{\e_n, 1} J_{\e_n}' [{\bf P}_{\e_n}] + \mu_{\e_n, 2} \Phi_{\e_n}'[{\bf P}_{\e_n}]=0.$$

Last expression can be read as $J_{\e_n}' [{\bf P}_{\e_n}]$ and
$\Phi_{\e_n}'[{\bf P}_{\e_n}]$ are linearly dependent. Observe
that this contradicts either the smoothness of $\partial
\mathcal{D}_{\e_n}$ or the nondegeneracy of $J_{\e_n}' [{\bf
P}_{\e_n}]$ on the tangent space.

For the sake of clarity, in what follows we will drop the
subscript $n$. Moreover, at many steps of the arguments we will
pass to a subsequence, without further notice. By using
Proposition \ref{energy111} and \eqref{crucial}, we have:
\begin{equation} \label{sist1} \frac{c_2}{2}\sum_{i=1}^\ell \overline{M}[P_i^\e]^2+ \frac{c_3}{2}\sum_{i\neq
j}w\Big(\frac{P_i^\e-P_j^\e}{\e}\Big) = \frac{c_4}{2} \e^{2 \beta},
\end{equation}
\begin{equation} \label{sist2} \frac{c_2}{2}\sum_{i=1}^\ell M[P_i^\e]^2- \frac{c_3}{2}
\sum_{i\neq j}\tau_{i}\tau_{j}
w\Big(\frac{P_i^\e-P_j^\e}{\e}\Big)=o(\e^{2\beta})
\end{equation}
\begin{equation} \label{sist3} \begin{array}{c} \dis c_2 \Big( (\mu_{\e, 2}+\mu_{\e, 1})M^+[P_i^\e] +
(\mu_{\e, 2}-\mu_{\e, 1}) M^-[P_i^\e] \Big)\\ \\ +\dis \mu_{\e, 2}\frac{c_3}{\e}
\sum_{j,\,j\neq i}
w'\Big(\frac{|P_i^\e-P_j^\e|}{\e}\Big)\frac{P_{i}^\e-P_{j}^\e}{|P_i^\e-P_j^\e|}
 -  \mu_{\e, 1} \sum_{j,\,j\neq i} \frac{\tau_{i}\tau_{j}}{\e^N}
\frac{\partial}{\partial P_i} \int_{\R^N} f(w_{P_i^\e})w_{P_j^\e} dx=
o(\e^{\beta}), \ i = 1, \dots, \ell. \end{array}
\end{equation}

Combining \eqref{sist1} and \eqref{sist2} we obtain:
\begin{equation} \label{sist+} c_2\sum_{i=1}^\ell M^+[P_i^\e]^2+ c_3\sum_{\tau_i=-\tau_j}
w\Big(\frac{P_i^\e-P_j^\e}{\e}\Big) = \frac{c_4}{2} \e^{2 \beta} +
o(\e^{2\beta}),
\end{equation}
\begin{equation} \label{sist-} c_2\sum_{i=1}^\ell M^-[P_i^\e]^2+ c_3\sum_{\tau_i=\tau_j}
w\Big(\frac{P_i^\e-P_j^\e}{\e}\Big) = \frac{c_4}{2} \e^{2 \beta} +
o(\e^{2\beta}).
\end{equation}

Motivated by \eqref{sist3}, we distinguish two cases:

{\bf Case 1:} There exists $C>0$ independent of $\e$ such that
$\sum_{i=1}^{\ell} |\mu_{\e, 2}+\mu_{\e, 1}|M^+[P_i^\e]^2 + |\mu_{\e,
2}-\mu_{\e, 1}| M^-[P_i^\e]^2 \geq C\e^{2\beta}$.

For instance, we can assume that $\sum_{i=1}^{\ell} |\mu_{\e,
2}+\mu_{\e, 1}|M^+[P_i^\e]^2 \geq C\e^{2\beta}$. In particular, recalling that $M^+[P_i^\e]^2<\e^{2\beta}$,
this implies that $|\mu_{\e, 2}+\mu_{\e, 1}| \nrightarrow 0$ and there exists $i_{0}
\in \{ 1, \dots ,\ell \}$ such that $M^+[P_{i_{0}}^\e]^2 \geq C
\e^{2\beta}$.

The idea is the following: we make the derivative $\mu_{\e, 1}
J_{\e}' [{\bf P}_{\e}] + \mu_{\e, 2} \Phi_{\e}'[{\bf P}_{\e}]$
along the same direction for all points ``close" to $P_{i_{0}}^\e$.
Since the direction is the same, the derivative of the interaction among those points
should be zero. And this direction will be chosen
conveniently to get a contradiction.

Take $\beta' \in (\beta, 1)$ fixed; let us define
$$ I = \{i=1, \dots ,\ell:\ |P_i^\e-P_{i_0}^\e|=o(\e^{\beta'})\}.$$

We take ${P_{i_{0}}^\e}^{\!+}=(P_{i_0,1}^\e, \dots P_{i_0,r}^\e, 0, \dots, 0)$ the
projection of $P_{i_{0}}^\e$ onto $A$ (here $P_{i_0,n}^\e$ denotes the $n$-th component of $P_{i_0}^\e$). Recall that $|P_i^\e|
=O(\e^{\beta})$. By multiplying \eqref{sist3} by ${P_{i_{0}}^\e}^{\!+}$ and
adding in $i \in I$, we have:
\begin{equation} \label{facil}  \begin{array}{c} \dis \sum_{i \in I}  c_2
(\mu_{\e,2}+\mu_{\e,1})M^+[P_i^\e, P_{i_{0}}^\e]  \\ \\ +\dis
\sum_{i\in I}\sum_{j\neq i}\bigg[\mu_{\e,2} \frac{c_3}{\e}
w'\Big(\frac{|P_i^\e-P_j^\e|}{\e}\Big)\frac{P_{i}^\e-P_{j}^\e}{|P_i^\e-P_j^\e|} - \mu_{\e, 1} \frac{\tau_{i}\tau_{j}}{\e^N}
\frac{\partial}{\partial P_i} \int_{\R^N} f(w_{P_i^\e})w_{P_j^\e} dx
\bigg ]\cdot{P_{i_0}^\e}^{\!+}\!= o(\e^{2\beta}).
\end{array} \end{equation}

We now estimate each of the above terms in order to get a
contradiction. First, observe that $M^+[P_i^\e,P_{i_{0}}^\e] =
M^+[P_{i_{0}}^\e,P_{i_{0}}^\e] + M^+[P_i^\e-P_{i_{0}}^\e, P_{i_0}^\e] \geq C
\e^{2\beta} + o(\e^{\beta + \beta'})$. Recall also that $|\mu_{\e,
2}+\mu_{\e, 1}| \nrightarrow 0$. So, it suffices to show that the
rest of the terms in \eqref{facil} are negligible to obtain a
contradiction.

We split the second  sum in two terms; those with $j \in I$ and those with
$j \notin I$. Let us start with the latter; by using Lemma \ref{gamma0} we have:

\begin{equation} \label{negli1}\sum_{i\in I} \sum_{j\notin I} \bigg[\mu_{\e,2}\frac{c_3}{\e}
w'\Big(\frac{|P_i^\e-P_j^\e|}{\e}\Big)\frac{P_{i}^\e-P_{j}^\e}{|P_i^\e-P_j^\e|} - \mu_{\e, 1}\frac{c_3+o(1)}{\e} \tau_{i}\tau_{j}
w'\Big(\frac{|P_i^\e-P_j^\e|}{\e}\Big)\frac{P_{i}^\e-P_{j}^\e}{|P_i^\e-P_j^\e|}\bigg]\cdot
{P_{i_{0}}^\e}^{\!+}. \end{equation}

Observe that, by definition of $I$,
$|P_j^\e- P_i^\e| \geq C \e^{\beta'}$ for $i\in I,\ j \notin I$. This
implies that $w'\big(\frac{|P_i^\e-P_j^\e|}{\e}\big)=o(  e^{-C
\e^{-1+\beta'}})$, and then \eqref{negli1} is negligible.

We now consider the following sum in $j \in I$:
$$\sum_{i \in I} \sum_{j\in I,\, j\neq i} \bigg[\mu_{\e,2} c_3 \frac{\partial}{\partial
P_i}w\Big(\frac{|P_i-P_j|}{\e}\Big) -\mu_{\e,
1}  \frac{\tau_{i}\tau_{j} }{\e^N} \frac{\partial}{\partial
P_i} \int_{\R^N} f(w_{P_i})w_{P_j} dx\bigg].$$

By a change of variables we deduce that $\int_{\R^N}
f(w_{P_i})w_{P_j} dx$ is a function of $|P_i-P_j|$. And it is easy
to conclude that for any $\xi \in {\cal C}^1(\R)$,
$$ \sum_{i \in I} \sum_{j\in I ,\, j\neq i} \frac{\partial}{\partial P_i}
\xi (|P_i-P_j|) = \sum_{i \in I} \sum_{j\in I,\,
j\neq i} \xi'(|P_i-P_j|) \frac{P_i-P_j}{|P_i-P_j|}=0.
$$
So we get a contradiction in Case 1.

\bigskip

{\bf Case 2:} $\sum_{i=1}^{\ell} |\mu_{\e, 2}+\mu_{\e, 1}|
M^+[P_i^\e]^2 + |\mu_{\e, 2}-\mu_{\e, 1}| M^-[P_i^\e]^2 =
o(\e^{2\beta})$. \bigskip

In a sense, here the effect of $M[P_i^\e]^2$ is negligible and the
interaction among the bumps is important. But in \cite{dapi2} it
was proved that the bumps cannot reach an equilibrium by
themselves (see Lemma \ref{lema42}), and this gives us the desired
contradiction.

Since $\mu_{\e, 1}^2+\mu_{\e, 2}^2=1$, then at least one between $ \mu_{\e, 1}+\mu_{\e, 2}$ and $\mu_{\e, 1}-\mu_{\e, 2}$ does not go to $0$. If $\mu_{\e, 2}+\mu_{\e, 1} \not\to 0 \Rightarrow \sum_{i}
M^+[P_i^\e]^2 = o(\e^{2\beta})$, and by \eqref{sist+},
$\sum_{\tau_i=-\tau_j} w\big(\frac{P_i^\e-P_j^\e}{\e}\big) \geq C \e^{2\beta}$. Analogously, if
$\mu_{\e, 2}-\mu_{\e, 1} \nrightarrow 0$ we can use \eqref{sist-}
to conclude $\sum_{\tau_i=\tau_j} w\big(\frac{P_i^\e-P_j^\e}{\e}\big) \geq C\e^{2\beta}$.

In any case, we have:
$$ \sum_{i \neq j} |\mu_{\e,2} - \mu_{\e,1} \tau_i \tau_j| w\Big(\frac{P_i^\e-P_j^\e}{\e}\Big)\geq C\e^{2\beta} .$$ So, there exist  $i_0 \neq j_0$ so that
$$|\mu_{\e,2} - \mu_{\e,1} \tau_{i_0} \tau_{j_0}|
w\Big(\frac{P_{i_0}^\e-P_{j_0}^\e}{\e}\Big) \geq C \e^{2\beta} .$$

Define:
$$ I = \left\{i=1, \dots, \ell:\ \frac{|P_i^\e-P_{i_0}^\e|}{\e \log (1/\e)} \mbox{ is bounded}\right\}.$$

Observe that, at least, $i_0,\ j_0 \in I$. For any $i \in I$, we
can pass to the limit on the following expressions:
$$ \frac{P_i^\e-P_{i_0}^\e}{2 \beta \e \log (1/\e)} \longrightarrow  Q_i \in \R^N, $$ and
$$ \e^{-2\beta} (\mu_{\e,2} - \mu_{\e,1} \tau_i \tau_j)
w\Big(\frac{P_i^\e-P_j^\e}{\e}\Big) \longrightarrow a_{i,j} \in \R,\; \ j = 1 ,\dots,
\ell.$$

We point out that $Q_{i_0}=0$ and $a_{i_0, j_0} \neq 0$. We recall that $w\big(\frac{P_i^\e-P_j^\e}{\e}\big) = O(\e^{2\beta})$; hence, from
\eqref{wdecay}, we obtain: \begin{equation} \label{Qi} Q_i -
Q_j = \lim_{\e \to 0} \frac{P_i^\e-P_{j}^\e}{2 \beta \e \log (1/\e)}
\Longrightarrow |Q_i - Q_j| \geq 1,\; \ i, \ j \in I.\end{equation}

Before going on, we are interested in extracting consequences from
$a_{i,j}\neq 0$, with $i \in I$. In such case there exist $c$,
$c'$ positive constants such that $c \, \e^{2\beta} \geq
w\big(\frac{P_i^\e-P_j^\e}{\e}\big)\geq c' \e^{2\beta}$. Then,
$$\frac{|P_i^\e-P_{j}^\e|}{2 \beta \e \log (1/\e)} \longrightarrow 1.$$
In particular, $j \in I$. Moreover, similarly as in \eqref{Qi}, we
obtain that $ |Q_i - Q_j| = 1$.

By using \eqref{sist3} together with Lemma \ref{gamma0}, we get: $$ \begin{array}{c} \dis c_2 \Big(
(\mu_{\e, 2}+\mu_{\e, 1})M^+[P_i^\e] + (\mu_{\e, 2}-\mu_{\e, 1})
M^-[P_i^\e] \Big)+ \\ \\ \dis \frac{c_3}{\e} \sum_{j,\,j\neq i} \big(\mu_{\e, 2}-\mu_{\e, 1} (1+o(1))\tau_i
\tau_j\big)
w'\Big(\frac{|P_i^\e-P_j^\e|}{\e}\Big)\frac{P_{i}^\e-P_{j}^\e}{|P_i^\e-P_j^\e|}=o(\e^{\beta}),\;\; \ i \in I. \end{array}
$$
We multiply by $\e^{1-2\beta}$ and use \eqref{wdecay} to obtain

$$ \dis  c_3 \sum_{j,\,j\neq i} \e^{-2\beta}(\mu_{\e, 2}-\mu_{\e, 1} \tau_i
\tau_j)
w\Big(\frac{P_i^\e-P_j^\e}{\e}\Big)\frac{P_{i}^\e-P_{j}^\e}{|P_i^\e-P_j^\e|} =
o(1), \;\;\ i \in I.$$

Recall that $a_{i,j}=0$ for any $j \notin I$. Passing to the
limit:
\begin{equation} \label{d-crit}\sum_{j\in I,\,j\neq i} a_{i,j}
\frac{Q_{i}-Q_{j}}{|Q_i-Q_j|}=0, \ i \in I.\end{equation}

In other words, the points $Q_i \in \R^N$, $i \in I$, satisfy that
$Q_{i_0}=0$, $|Q_i-Q_j| \geq 1$ and $(Q_i)_{i\in I}$ is a critical
point of the function:
\begin{equation}\label{crit}  (Z_i)_{i\in I} \longmapsto \sum_{i,j\in I,\,i\neq j} a_{i,j} |Z_i-Z_j|,\quad Z_i\in\R^N. \end{equation} where $a_{i,j}=a_{j,i}$,  $a_{i,j} = 0$
for points $Q_i$, $Q_j$ such that $|Q_i-Q_j|>1$, and $a_{i_0, j_0}
\neq 0$.

We finish the proof by showing that this is impossible. For that
we need to distinguish between the case of positive peaks and the
case of mixed positive and negative peaks.

In the first case $\tau_i =1 $ for all $i=1,\ldots, \ell$. By the definition of
$a_{i,j}$ and the fact $a_{i_0,j_0} \neq 0$, we conclude that
$\mu_{\e,2} - \mu_{\e, 1} \nrightarrow 0$. Moreover, $a_{i,j}$
have all the same sign as $\mu_{\e,2} - \mu_{\e, 1}$. Assume, for
instance, that $a_{i,j} \geq 0$. But in such case $(Q_i)_{i\in I}$
cannot be a critical point of the map given by \eqref{crit}, as
can be seen using dilatations. More specifically, if we multiply
\eqref{d-crit} by $Q_i$ and make the addition, we get:
$$  \sum_{i \in I} \sum_{j\in I,\,j\neq i} a_{i,j}
\frac{Q_{i}-Q_{j}}{|Q_i-Q_j|}\cdot Q_i = \sum_{i<j} a_{i,j}
|Q_{i}-Q_{j}| \geq a_{i_0,j_0}>0. $$

The case in which there are peaks of different sign is excluded
thanks to the next lemma, proved in \cite{dapi2}. We point out
that the restriction $\ell \leq 6$ is needed only at this point.

\begin{lemma} \label{lema42} Let $\ell \geq 2$ and consider the
function:
$$ \Phi: (Z_1, \dots, Z_{\ell}) \in \R^{N \ell} \to \sum_{i \neq j} a_{ij} |Z_i - Z_j|.  $$ where $a_{ij}=a_{ji}$. Suppose that $\Phi$ is not
identically zero and that there exists  a critical point $({Q}_1, \dots,
{Q}_{\ell})$ of $ \Phi$ satisfying:
$$ |Q_i-Q_j|\geq 1 \mbox{ for }i\neq j\;\;\;\mbox{ and }\;\; \;|Q_i-Q_j|=1 \mbox{ if } a_{i,j}\neq 0. $$
Then $\ell \geq 7$.

\end{lemma}

\noindent{\bf Proof of Theorems \ref{th1} and \ref{th2} completed.}  According to Lemma \ref{relation}, for
$\e>0$ sufficiently small $\chi w_{{{\bf P}}_\e}+ \phi_{{{\bf P}}_\e}$ solves the equation \eqref{sch},
 where ${{\bf P}}_\e=(P_1^\e,\ldots, P_\ell^\e)\in \Gamma_\e$ is the critical
 point of $J_\e$ with critical value $J^*_\e$.
 The construction of the family ${\bf P}_\e$
 depends on the particular $\beta\in (0,1)$ chosen at the beginning of Section 2.
 To emphasize this fact we denote this family as ${{\bf P}}_{\e,\beta}$.
 Let $\beta_k\subset (0,1)$ be any sequence such that $\beta_k\to 1$. Then there is a
 decreasing sequence of positive numbers $\e_k$ such that for all $0<\e<\e_k$ one has:

 \begin{enumerate}

 \item $\chi w_{{{\bf P}}_{\e,\beta_k}}+ \phi_{{{\bf P}}_{\e,\beta_k}}$ solves
  \eqref{sch},

 \item by \eqref{GMC}, $|P_i^{\e,\beta_k}|\leq (\min_{i} |\lambda_i| )^{-1/2} \, \e^{\beta_k}$,
  $|{{P}}_i^{\e,\beta_k}-{{P}}_j^{\e,{\beta_k}}|\geq
  2\beta_k^2\e\log\frac{1}{\e}$ for $i\neq j$, and

  \item $|\phi_{{{\bf P}}_{\e,\beta_k}}|\leq
  \e^{\beta_k^2(1+\sigma)}$.

\end{enumerate} We define ${{\bf P}}_\e={{\bf P}}_{\e,\beta_k}$ and $v_\e=\chi
   w_{{{\bf P}}_{\e,\beta_k}}+ \phi_{{{\bf P}}_{\e,\beta_k}}$  if $\e_{k+1}<\e<\e_k$
   and we clearly have that the theses of Theorem \ref{th1} and Theorem \ref{th2} hold.

\appendix
\renewcommand{\theequation}{\Alph{section}.\arabic{equation}}
\section{Key energy estimate}
Consider the configuration set $\Gamma_\e$ and the approximate
solutions $\chi w_{\bf P}$ defined in Section 2. In this Appendix
we will derive some crucial estimates. We note that by assumption
(V2) and  \eqref{GMC}  we have $|\nabla V(P_i)|\leq C\e^{\beta}$
for
 ${\bf P}\in \Gamma_\e$; then  by \eqref{wdecay}
 we deduce
$$\begin{aligned}|V(x)\chi w_\epii-V(P_i) w_\epii|&\leq|\nabla V(P_i)||x-P_i|w_\epii+ C|x-P_i|^2w_\epii\leq C\e^{1+\beta}w_\epii^{2/3},
\end{aligned}$$
by which
\begin{equation}\label{andals}V(x)\chi w_\epii-V(P_i)w_\epii=O(\e^{1+\beta})w^{2/3}_\epii,\qquad V(x)\chi
w_\pii- w_\pii=O(\e^{2\beta})w_\pii^{2/3}\end{equation} uniformly
for
 ${\bf P}\in \Gamma_\e$.

\begin{remark}\label{gamma00}
 Observe that by \eqref{wdecay} it follows that
$$\frac{w(z+\xi)}{w(\xi)}\leq C \bigg(\frac{|\xi|}{1+|z+\xi|}\bigg)^{\frac{N-1}{2}}e^{|z|}
\quad \forall z,\xi\in \R^N.$$ Since $\frac{|\xi|}{1+|z+\xi|}\leq
2(1+|z|)$, we deduce
\begin{equation}\label{domdecay}w(z)\frac{w(z+\xi)}{w(\xi)}\leq C\quad \forall z,\xi\in \R^N.\end{equation}
By taking $z=\frac{x-P_i}{\e}$ and $\xi=\frac{P_i-P_j}{\e}$, \eqref{domdecay} yields
\begin{equation}\label{esti00}w_\epii w_\epij
\leq Cw\Big(\frac{P_i-P_j}{\e}\Big) \leq  C \e^{2\beta}
\quad\forall i\neq j\end{equation} uniformly for ${\bf P}\in
\Gamma_\e$.
\end{remark}

The next two lemmas are devoted  to estimate some integrals
associated to $w_\epii$'s.

\begin{lemma} \label{gamma0}
For $i\neq j$ the following expansions hold uniformly for ${\bf
P}\in {\Gamma}_\e$:
\[
 \intr
  f(w_\epii)
w_\epij\,dx=c_3\e^N(1+o(1))w\Big(\frac{P_i-P_j}{\e}\Big),
\]
\[
 \frac{\partial}{\partial P_i}\bigg[\intr
  f(w_\epii)
w_\epij\,dx\bigg]=c_3\e^{N-1}(1+o(1))w'\Big(\frac{|P_i-P_j|}{\e}\Big)\frac{P_i-P_j}{|P_i-P_j|}.
\]
where $c_3=\intr f(w)e^{
x_1} dx.$
\end{lemma}
 \noindent{\textit{Proof.}} First
consider the function
$$\xi(\rho)=\intr f(w)w(x+\rho {\bf e}_1 )dx,\quad\rho>0,$$ where
${\bf e}_1$ is the first vector of the standard basis of $\R^N$,
i.e. ${\bf e}_1=(1,0,\ldots,0)$.
 According to
(\ref{wdecay}) for every $x\in\R^N$ we have
\begin{equation}\label{doco}\lim_{\rho\rightarrow \infty}\frac{w(x+\rho {\bf e}_1)}{w(\rho)}=\lim_{
\rho\to\infty}e^{-|x+\rho {\bf e}_1|+\rho}=
e^{-x_1}.\end{equation} Thanks to  \eqref{domdecay} the  Dominated
Convergence Theorem applies and gives
$\frac{\xi(\rho)}{w(\rho)}\to \int_{\R^N}f(w)e^{-x_1} dx $. Next
compute
$$ \xi'(\rho)=\intr f(w)w'(x+\rho {\bf e}_1 )\frac{x_1+\rho}{|x+\rho {\bf e}_1|}
dx.$$ Using \eqref{wdecay} and proceeding as above we get
$$\frac{\xi'(\rho)}{ w'(\rho)}\to\intr f(w)e^{-x_1}dx.$$
Since
$$\intr f(w_\epii)w_\epij dx=\e^N\intr f(w)w\Big(x+\frac{P_i-P_j}{\e}\Big) dx=\e^N\xi\Big(\frac{|P_i-P_j|}{\e}\Big),$$
and
$$\frac{\partial }{\partial P_{i}}\bigg[\intr f(w_\epii)w_\epij dx\bigg]=\e^{N-1}\xi'\Big(\frac{|P_i-P_j|}{\e}\Big)\frac{P_i-P_j}{|P_i-P_j|},$$
then the thesis follows. \hfill $\Box$

\bigskip

\begin{lemma}\label{gamma000}
For every $i=1,\ldots, \ell$  the following asymptotic expansion
holds uniformly for ${ \bf P}\in{\Gamma}_\e$:
$$\intr V(x) \chi^2 w_\epii \nabla w_\epii dx=-\frac{\e^N}{2}M[P_i]\intr w^2 dx  +o(\e^{N+\beta}).$$
\end{lemma}
\noindent {\textit{Proof.}}
Observe that
$$|\nabla(\chi^2
V(x))-\nabla V(P_i )-D^2 V(P_i )(x-P_i )|w_{\epii }\leq
C|x-P_i|^2w_{\epii }\leq C\e^2w_{\epii }^{1/2}$$ by which, using
that  $\intr  D^2 V(P_i )(x-P_i )  w_{\epii }^2 dx =\e^N\intr D^2
V(P_i )y  w^2(y) dy=0$,

\begin{equation*}\begin{aligned}
2\intr  V(x) \chi^2 w_\epii \nabla w_\epii dx&=-\intr \nabla
(\chi^2 V)(x)w_\epii^2dx=-\intr \nabla
V(P_i)w_\epii^2dx+O(\e^{N+2})\\ &=-\e^N\nabla V(P_i)\intr w^2 dx
+O(\e^{N+2})
\end{aligned}\end{equation*}
uniformly for ${ \bf P}\in{\Gamma}_\e$, and, since $ \nabla V(P)
=M[P]+o(|P|)$ as $P\to 0$, we obtain the thesis. \hfill $\Box$

\bigskip

 The next proposition provides an estimate of the error
up to which the functions $\chi w_\epi$ satisfy \eqref{sch}.
\begin{lemma}\label{error}
There exists a constant $C>0$ such that for every $\e>0$ and ${\bf
P}\in \Gamma_\e$:
\begin{equation*}
 \big|{\cal S}_\e[\chi w_\epi] |
\leq C\e^{\beta(\beta+\sigma)}\sum_{i=1}^{\ell}w_{\epii}^{1-\beta}
\end{equation*} where ${\cal S}_\e$ is the operator defined in \eqref{operator}.
\end{lemma}

\noindent {\textit{Proof.}} By \eqref{andals}  we deduce
\begin{equation*} \begin{aligned}\e^2\Delta (\chi w_{ {\bf P}}) - V(x)\chi w_{{\bf P}}
+ f(\chi w_{{\bf P}}) &= \e^2\Delta w_{ {\bf P}} - w_{{\bf P}} +
f(w_{{\bf P}}) +O(\e^{2\beta}) \sum_{i=1}^{\ell}w_{\epii}^{2/3}\\
&=f(w_{{\bf P}})-\sum_{i=1}^{\ell}
\tau_{i}f(w_{P_i})+O(\e^{2\beta})
\sum_{i=1}^{\ell}w_{\epii}^{2/3}\end{aligned}\end{equation*}
uniformly for  ${\bf P }\in\Gamma_\e$. Given ${\bf
P}\in\Gamma_\e$, in the following we will make use of the
following sets $A_{\e,i}$
$$A_{\e,i}=\Big\{x\in\R^N\,\Big|\,w_\epii>a\e^{\beta}\Big\}$$ where $a>0$ is chosen such that, according to \eqref{esti00},
$$A_{\e,i}\cap A_{\e,j}=\emptyset \;\;\forall i\neq j.$$
Observe that
\begin{equation}\label{kapp1}w_\epij\leq a \e^{\beta}\leq w_\epii
\hbox{ on }A_{\e,i}\hbox{ for }j\neq i.\end{equation} Then, by
using assumption (f1), we get \begin{equation*}\Big|f(w_{\bf
P})-\tau_{i} f(w_\epii)\Big|\leq C w_{P_i}^{\sigma}\sum_{j\neq i}
w_{\epij} \;\hbox{ on }A_{\e,i},\end{equation*} by which
$$\Big|f(w_{\bf
P})-\tau_{i} f(w_\epii)\Big|\leq
C\e^{\beta(\beta-\sigma)}\sum_{j\neq i} (w_\epii
w_\epij)^{\sigma}w_\epij^{1-\beta}\leq
C\e^{\beta(\beta+\sigma)}\sum_{j\neq i}w_j^{1-\beta} \;\hbox{ on
}A_{\e,i}.$$ On the other hand
\begin{equation*}\begin{aligned}|f(w_\epii)| \leq C |w_\epii|^{1+\sigma}\leq C \e^{\beta(\beta+\sigma)}w_\epii^{1-\beta}\hbox{ on }\R^N\setminus  A_{\e,i},
\end{aligned}\end{equation*}
\begin{equation*}|f(w_\epi)|\leq C\sum_{j=1}^\ell |w_\epij|^{1+\sigma}\leq C \e^{\beta(\beta+\sigma)}\sum_{j=1}^\ell w_\epij^{1-\beta}\hbox{ on }\R^N\setminus \cup_{j=1}^\ell A_{\e,j}.\end{equation*}
Since $\beta(\beta+\sigma)<2\beta$ we obtain the thesis.
\hfill$\Box$

\bigskip

 With the help of Lemma \ref{gamma0} and Lemma \ref{gamma000} we derive the following  key energy estimate.

 \begin{proposition}
The following asymptotic expansions  hold uniformly for ${\bf P
}=(P_1,\ldots,P_\ell)\in {\Gamma}_\e$:
\begin{equation}\label{lonel}\begin{aligned}I_\e [\chi w_{{\bf
P}}]=& c_1\e^N
+\frac{c_2}{2}\e^N\sum_{i=1}^{\ell}M[P_i]^2-\frac12 \sum_{i\neq
j}\tau_i\tau_j\intr f(w_\epii )w_\epij dx +o(\e^{N+2\beta})
,\end{aligned}
\end{equation}
$$\frac{\partial }{\partial P_i}\Big(I_\e [\chi w_{{\bf
P}}]\Big)=c_2\e^N M[P_i]-\sum_{j,\, j\neq
i}\tau_{i}\tau_{j}\frac{\partial }{\partial P_i}\bigg(\intr f(w_\epii )w_\epij dx\bigg)
+o(\e^{N+\beta}),\quad
i=1,\ldots, \ell,$$
 where the constants $c_1,c_2$ are given by
$$c_1=\frac{\ell}{2}\intr\big( |\nabla w|^2+w^2\big)dx-\ell\intr F(w)dx,\quad c_2=\frac12\intr w^2dx. $$ \label{energy1}
\end{proposition}
\noindent{\textit{Proof.}} We begin by estimating the potential
term: by  \eqref{andals} we derive
\begin{equation}\label{litti1}\begin{aligned}  \intr  V(x)|\chi w_\epi|^2 dx &=\sum_{i=1}^\ell\intr V(x) |\chi w_\epii|^2
  dx+\sum_{i\neq  j}\tau_{i}\tau_j\intr V(x) \chi w_\epii \chi w_\epij dx\\ &
 =\sum_{i=1}^\ell V(P_i)\e^N \intr w^2
  dx+\sum_{i\neq j}\tau_{i}\tau_jV(P_i) \intr  w_\epii
  w_\epij dx+o\big(\e^{N+2\beta}\big)\\ &=\sum_{i=1}^\ell \Big(1+\frac12M[P_i]^2\Big)\e^N \intr w^2
  dx+\sum_{i\neq j}\tau_{i}\tau_j \intr  w_\epii
  w_\epij dx+o\big(\e^{N+2\beta}\big)
\end{aligned}\end{equation}uniformly for ${ \bf P}\in{\Gamma}_\e$, where the last equality follows by assumption (V2) and \eqref{esti00}.

Next we compute
\begin{equation}\label{pis2}\begin{aligned}\frac{\e^2}{2}&\intr
|\nabla (\chi w_\epi)|^2 dx-\intr F(\chi w_\epi) dx
=\frac{\e^2}{2}\intr |\nabla w_\epi|^2
dx-\intr F( w_\epi) dx+o(\e^{N+2})\\
&=\ell\frac{\e^N}{2}\intr |\nabla w|^2 dx -\ell\e^N\intr F(w) dx\\
& \;\;\;\;+\frac{\e^2}{2}\sum_{i\neq j}\tau_{i}\tau_{j}\intr
\nabla w_\epii\nabla w_\epij dx-\intr \Big(F(
w_\epi)-\sum_{j=1}^\ell F(w_\epij)\Big)
dx+o(\e^{N+2})\end{aligned}
\end{equation}uniformly for ${\bf
P}\in{\Gamma}_\e$.

 Combining  \eqref{litti1} with \eqref{pis2}, and using equation \eqref{ground},  we
get
\begin{equation}\label{pd6}\begin{aligned} I_{\e}[\chi w_{{\bf
P}}] &= c_1\e^N +\frac{c_2}{2} \e^N\sum_{i=1}^\ell M[P_i]^2 +\fr\sum_{i\neq
j}\tau_{i}\tau_j
\intr f(w_\epii)w_\epij dx \\
& \;\;\;\;-\intr \Big(F( w_\epi)-\sum_{i=1}^\ell F(w_\epii)\Big)
dx+o(\e^{N+2\beta})\\&=c_1\e^N +\frac{c_2}{2} \e^N\sum_{i=1}^\ell
M[P_i]^2 -\fr\sum_{i\neq j}\tau_{i}\tau_j \intr f(w_\epii)w_\epij
dx -H({\bf P})+o(\e^{N+2\beta}),\end{aligned}
\end{equation}
uniformly for ${\bf P}\in{\Gamma}_\e$, where we have set
$$H({\bf P})=\intr  F(w_{\bf P}) dx
-\sum_{i=1}^\ell\intr F(w_{\epii})dx -\sum_{i\neq j}\tau_{i}\tau_j
\intr f(w_{\epii}) w_{\epij} dx ,\quad {\bf P}\in \Gamma_\e.$$
Consider the sets $A_{\e,i}$ defined in Lemma \ref{error}; by
assumption (f1) we have
\begin{equation*}\begin{aligned} |H({\bf P})|&\leq \sum_{i=1}^\ell
\int_{A_{\e,i}}\!\bigg|F(w_{\bf P}) -F(w_{\epii})
-f(w_\epii)\sum_{j\neq i}\tau_{i}\tau_{j}w_\epij \bigg| dx  +C\!\sum
_{i=1}^\ell \int_{\R^N\setminus A_{\e,i}}\!\!
\bigg(w_\epii^{2+\sigma}\!+w_\epii^{1+\sigma}\sum_{j\neq i}w_\epij
\bigg)dx .\end{aligned}
\end{equation*}
By \eqref{kapp1} we get \begin{equation*}\Big|F(w_{\bf
P})-F(w_\epii)-f(w_\epii)\sum_{j\neq
i}\tau_{i}\tau_{j}w_{P_j}\Big|\leq C w_{P_i}^\sigma\sum_{j\neq i}
w_{\epij}^{2}=C \sum_{j\neq i} (w_{P_i}w_\epij)^\sigma
w_{\epij}^{2-\sigma} \;\hbox{ on }A_{\e,i}.\end{equation*} Taking
into account of  \eqref{esti00} and \eqref{kapp1}, the above
inequalities  imply $H({\bf P})=o(\e^{N+2})$ uniformly for ${\bf
P}\in\Gamma_\e$. Then by  \eqref{pd6}
 we obtain \eqref{lonel}.

We now estimate  the error  term $o(\e^{N+2\beta})$ in
\eqref{lonel} in the $\mathcal C^1$ sense. To this aim, fix
$i\in\{1,\ldots, \ell\}$;  by definition we have  $\frac{\partial
w_{\epi}}{\partial P_{i}}=-\tau_{i} \chi\nabla w_{\epii}$.
 Then, by using Lemma \ref{gamma000}, we can compute
\begin{equation*}\begin{aligned}\tau_{i}\frac{\partial I_\e[\chi
w_\epi]}{\partial P_{i}}&=-\big\langle I'_\e[\chi
w_\epi],\chi\nabla w_\epii\big\rangle=\intr \big(\e^2\Delta(\chi
w_\epi)-V(x)\chi w_\epi+f(\chi w_\epi)\big)\chi\nabla w_\epii dx
\\&=c_2\tau_i \e^NM[P_i]+\intr\Big(\e^2\sum_{ j=1}^\ell\tau_{j}\Delta w_\epij-\sum_{ j\neq i} \tau_{j}\chi V(x) w_\epij+f(w_\epi)\Big)\nabla w_\epii dx+o(\e^{N+\beta})
\end{aligned}\end{equation*}
uniformly for ${\bf P}\in{\Gamma}_\e$. Since $|\nabla w_\epii|\leq
C \e^{-1}w_\epii$, using \eqref{andals} and \eqref{esti00}, for
$i\neq j$ we have
$$\intr \chi
V(x) w_{\epij }\nabla w_{\epii }\,dx=\intr w_{\epij }\nabla
w_{\epii }+ o(\e^{N+\beta}),$$ while $\intr w_\epii\nabla w_\epii
dx=\frac12\intr \nabla w_\epii^2dx =0$. Then, using \eqref{ground}
we arrive to
\begin{equation*}\begin{aligned}\tau_{i}\frac{\partial I_\e[\chi
w_\epi]}{\partial P_i}
&=c_2\tau_i\e^NM[P_i]+\intr\Big(\e^2\sum_{j=1}^\ell\tau_{j} \Delta
w_\epij-\sum_{j=1}^\ell \tau_{j}w_\epij+f(w_\epi)\Big)\nabla
w_\epii dx+o(\e^{N+\beta})\\
&=c_2\tau_i\e^NM[P_i]+\intr\Big(f(w_\epi)-\sum_{j=1}^\ell\tau_{j}f(w_\epij)
\Big)\nabla w_\epii dx+o(\e^{N+\beta})
\\ &=c_2\tau_i\e^NM[P_i]+
\sum_{j\neq i}\tau_{j}\intr f'(w_\epii)w_\epij\nabla w_\epii dx +
K({\bf P})+o(\e^{N+\beta})\end{aligned}\end{equation*} uniformly
for ${\bf P}\in{\Gamma}_\e$, where we have set $$K({\bf
P})=\intr\Big(f(w_\epi)-\sum_{j=1}^\ell\tau_{j}f(w_\epij)-f'(w_\epii)\sum_{j\neq
i} \tau_{j}w_\epij\Big)\nabla w_\epii dx .$$
 By assumption (f1) it follows
\begin{equation}\label{kappa}\begin{aligned}|K({\bf P})|\leq& \int_{A_{\e,i}}\bigg|f(w_{\bf P}) -\tau_{i}f(w_{\epii})
-f'(w_\epii)\sum_{j\neq i}\tau_{j}w_\epij \bigg|\big|\nabla
w_\epii\big| dx\\ &+\sum_{j\neq i}\int_{A_{\e,j}}\big|f(w_{\bf P})
-\tau_{j}f(w_{\epij})\big|\big|\nabla w_\epii\big|\, dx\\ &+
C\sum_{j=1}^\ell\int_{\R^N\setminus A_{\e,j}}
w_\epij^{1+\sigma}\big|\nabla w_\epii\big| +C\sum_{j\neq
i}\int_{\R^N\setminus A_{\e,i}}w_\epii^\sigma w_\epij\big|\nabla
w_\epii\big|dx.
\end{aligned}\end{equation}
Observe that  \begin{equation}\label{kappa1}\Big|f(w_{\bf
P})-\tau_{i}f(w_\epii)-f'(w_\epii)\sum_{j\neq
i}\tau_{j}w_{P_j}\Big|\big|\nabla w_\epii\big|\leq
C\e^{-1}\sum_{j\neq i}
w_{\epij}^{1+\sigma}w_\epii=C\e^{-1}\sum_{j\neq i} (w_\epii
w_{\epij})w_\epij^\sigma.\end{equation} Next fix $j\neq i$: by
\eqref{kapp1} we have
\begin{equation}\label{kappa2}|f(w_\epi)-\tau_{j}f(w_\epij)|\big|\nabla w_\epii\big|\leq C \e^{-1}w_\epij^\sigma w_\epii\sum_{k\neq j}w_\epik=C \e^{-1}(w_\epii w_\epij)^\sigma w_\epii^{1-\sigma}\sum_{k\neq j}w_\epik \;\hbox{ on }A_{\e,j}.\end{equation}
Inserting \eqref{kappa1}-\eqref{kappa2} into \eqref{kappa}, and
using \eqref{esti00} and \eqref{kapp1}, we deduce  $K({\bf
P})=o(\e^{N+1})$ uniformly for ${\bf P}\in{\Gamma}_\e$. Thus we
have obtained
\begin{equation*}
\begin{aligned}\frac{\partial I_\e[\chi
w_\epi]}{\partial P_{i}} &=c_2\e^NM[P_i]+\sum_{j\neq
i}\tau_{i}\tau_{j}\intr f'(w_\epii)w_\epij \nabla w_\epii dx +
o(\e^{N+\beta})\\ &=c_2\e^NM[P_i]-\sum_{j\neq
i}\tau_{i}\tau_{j}\frac{\partial }{\partial P_{i}}\intr
f(w_\epii)w_\epij dx+o(\e^{N+\beta})\end{aligned}\end{equation*}
uniformly for ${\bf P}\in{{\Gamma}}_\e$, and the second part of
the  thesis follows. \hfill $\Box$

\bigskip

\section{Lyapunov-Schmidt Reduction}
In this appendix we carry out the reduction procedure sketched in
Section 3. In particular we will prove Lemma \ref{reg} and
Proposition \ref{energy111}. A large part of the proofs follows in
a standard way but we include some details here for completeness.
\subsection{The linearized equation}
Consider the functions $Z_{P_i,n}$ defined in Section 2. Observe
that by proceeding as in the proof of \eqref{andals} we deduce
\begin{equation}\label{esti22}\begin{aligned}Z_{P_i,n}&=(1-\e^2\Delta)\frac{\partial w_\pii }{\partial
x_n}+O(\e^{2\beta-1})w_\pii ^{2/3}=f'(w_\pii ) \frac{\partial
w_\pii}{\partial x_n}+O(\e^{2\beta-1}) w_\pii ^{2/3}
\end{aligned}\end{equation}
uniformly for  ${\bf P }\in \Gamma_\e.$
 After integration by parts  it is immediate to prove that
   \begin{equation}\label{ort}\bigg(\phi, \frac{\partial(\chi w_{\pii})}{\partial
 x_n}\bigg)_\e = \intr \phi Z_{P_i,n}\, dx\quad \forall \phi\in H^1_V(\R^N),\end{equation}
 then orthogonality to the functions $\frac{\partial(\chi  w_\pii)}{\partial
 x_n}$ in $H^1_V(\R^N) $ with respect to the scalar product $(\cdot,\cdot)_\e$ is equivalent to orthogonality to $Z_{P_i,n}$ in $L^2(\R^N) $.
Hence we easily get\begin{equation}\label{dom}\intr
Z_{P_i,n}\frac{\partial (\chi w_{P_j})}{\partial x_{m}} \,dx=
\delta_{ij}\delta_{nm}\e^{N-2}\Big\|\frac{\partial w }{\partial
x_1}\Big\|_{H^1(\R^N)}^2+o(\e^{N-2}) .\end{equation}
   uniformly for ${\bf P }\in \Gamma_\e$ ($\delta_{ij}$
and $\delta_{nm}$  denoting the Kronecker's symbols), where
$\|v\|_{H^1(\R^N)}^2:=\intr (|\nabla v|^2+|v|^2) dx$.

Let $\mu\in (0,\sigma)$ be  a sufficiently small number  and
introduce the following weighted norm:
\begin{equation}
\label{starnorm} \|\phi \|_{*,{\bf P}}:= \sup_{x\in \R^N} \bigg(\sum_{i=1}^\ell
w_\epii(x)\bigg)^{-\mu}|\phi(x)|.
\end{equation}
We first consider  a  linear problem: given ${\bf P }\in
\Gamma_\e$ and  $ \theta \in L^2(\R^N)$, find a function
$\phi$ and constants $\alpha_{in}$ satisfying
\begin{equation}
\left\{\begin{aligned}& {\cal L}_{{\bf P }} [\phi]= \theta +
\sum_{i,n} \alpha_{in} Z_{P_i,n},
\\
&\phi \in H^2(\R^N)\cap H^1_V(\R^N) , \;\intr \phi Z_{P_i, n}
\,dx=0\hbox{ for } i=1,\ldots, \ell,\, n=1,\ldots, N,
\end{aligned}
\right. \label{linear}
\end{equation}
  where $$ {\cal L}_{{\bf P }} [\phi]:=\e^2 \Delta \phi-V(x)\phi + f' (\chi w_{{\bf P }}) \phi .$$

\begin{lemma}\label{apriorii} There exists a constant $C>0$ such that, provided that $\e$ is  sufficiently small, if ${\bf P }\in \Gamma_{\e}$ and   $(\phi, \theta,\alpha_{in})$ satisfies  (\ref{linear}), then
   $$|\alpha_{in}|\leq C(\e^{1+\sigma}\|\phi\|_{*, {\bf P}}+\e\|\theta\|_{*, {\bf P}}).$$

   \end{lemma}

   \noindent{\textit{Proof.}}
By multiplying the equation in \eqref{linear}
  by $\frac{\partial  (\chi w_{P_j})}{\partial x_{m}}$
  and integrating over $ \R^N$, we get
   \begin{equation}
   \label{ceqn}
    \sum_{i,n}\alpha_{in} \int_{\R^N} Z_{P_i,n}\frac{\partial (\chi w_{P_j})}{\partial x_{m}} dx=
    -\int_{\R^N}
    \theta
    \frac{\partial(\chi w_{P_j})}{\partial x_{m}}dx + \int_{\R^N}  {\cal L}_{{\bf P }} [\phi]
    \frac{\partial (\chi w_{P_j})}{\partial x_{m}}
    dx.
   \end{equation}
      \noindent First examine the left hand side of (\ref{ceqn}).
By using (\ref{dom})
   \begin{equation}\label{ex0}\bigg|\sum_{i,n}\alpha_{in}
   \int_{\R^N}
   Z_{P_i,n}\frac{\partial (\chi w_{P_j})}{\partial x_{m}} dx\bigg|\geq C\e^{N-2}|
   \alpha_{jm}|+o(\e^{N-2})\sum_{i,n}|\alpha_{i,n}|.\end{equation}

\noindent The first term on the right hand side of (\ref{ceqn})
can be estimated as
   \begin{equation}\label{ex1} \int_{\R^N}\Big| \theta \frac{\partial (\chi w_{P_j})}{\partial x_{m}} \Big|
   dx
    \leq C\| \theta \|_{*,{\bf P}}\int_{\R^N} |\nabla w_{P_j}| dx
\leq C\e^{N-1}\| \theta \|_{*,{\bf P}}.\end{equation} Finally, by
using   \eqref{esti22},
   \begin{align*} \bigg|\int_{\R^N}  {\cal L}_{{\bf P }} [\phi]  \frac{\partial (\chi w_{P_j})}{\partial
   x_{m}}dx\bigg|&=\bigg|\int_{\R^N}\phi\bigg[-Z_{P_j,m}+f'(\chi w_{{\bf P }})
\frac{\partial (\chi w_{P_j})}{\partial
   x_{m}}\bigg]\,dx\bigg|
  \\ &\leq C\|\phi\|_{*, {\bf P}} \intr\Big|(f'(w_{{\bf P }})-f'(w_{P_j}))
\frac{\partial  w_{P_j}}{\partial
   x_{m}}\Big|\,dx+C\e^{N+2\beta-2}\|\phi\|_{*, {\bf P}}
\\&\leq C\e^{-1}\|\phi\|_{*, {\bf P}}\sum_{i\neq j}
   \intr w_{\epii}^\sigma w_{\epij}\,dx+C\e^{N+2\beta-2}\|\phi\|_{*, {\bf P}}\\ &\leq C\|\phi\|_{*, {\bf P}}
   (\e^{N+2\beta\sigma-1}+\e^{N+2\beta-2}).
   \end{align*}
   where last inequality follows from \eqref{esti00}.
Combining this  with (\ref{ceqn}), (\ref{ex0}) and (\ref{ex1}), we
achieve the thesis. \hfill $\Box$

\bigskip

Now we prove  the following  a priori estimate for (\ref{linear}).

\begin{lemma}\label{apriori} There exists a  constant $C>0$ such that,
   provided that $\e$ is  sufficiently small, if ${\bf P }\in \Gamma_{\e}$ and   $(\phi, \theta,\alpha_{in})$ satisfies  (\ref{linear}), the following holds:
   $$\| \phi \|_{*,{\bf P}} \leq C \|\theta\|_{*,{\bf P}}
.$$
 \end{lemma}

   \noindent
   {\textit{Proof.}}  We argue by contradiction.   Assume the existence of a   sequence  $\e_k \to 0^+$,
   ${\bf P }_k \in
   \Gamma_{\e_k}$ and
   $ ({\phi}_{k},\theta_k, {\alpha}_{in}^k)$ satisfying
(\ref{linear})
   such that  \[\|\phi_{k}\|_{*,{\bf P}_k }=1, \;\; \|\theta_k\|_{*, {\bf P}_k }=o(1).\]
By Lemma \ref{apriorii} we deduce  $
   \alpha_{in}^k= o(\e)$ for every $(i,n)$, by which $\| \theta_k +\sum_{i,n} \alpha_{in}^kZ_{P_i^k ,n}\|_{*,{\bf P}_k } =o(1)$ and, consequently,
 \begin{equation} \label{del1}
 \| \e_k^2\Delta \phi_k -V(x)\phi_k +
  f' (\chi w_{{\bf P }_k}) \phi_k \|_{*,{\bf P}_k }=o(1).
\end{equation}
 We claim that \begin{equation}\label{cla} \| \phi_k
\|_{L^\infty (\cup_{i=1}^\ell B_{R\e_k} (P_i ^k))} = o(1)\quad
\forall R>0.\end{equation} Otherwise, we may assume that $
\|\phi_k\|_{L^\infty (B_{R\e_k} (P_1^k))} \geq c
>0$ for some $R>0$. By multiplying the equation in (\ref{linear}) by $\phi_{k}$ and
integrating by parts we immediately get that the sequence
$\phi_k(\e_k x+P_1^k)$ is bounded in  $H^1(\R^N)$. Therefore,
possibly passing to a subsequence,  $ \phi_k (\e_k x+P_1^k)
\rightharpoonup \phi_0 $ weakly in $H^1(\R^N)$ and a.e. in $\R^N$,
and $\phi_0$ satisfies
\begin{equation*} \Delta \phi_0-\phi_0 + f' (w) \phi_0=0,\;\;
|\phi_0(x)|\leq  w^{\mu}(x).
\end{equation*}
 According to elliptic regularity theory   we may assume $\phi_k(\e_k x+P_1^k)\to \phi_0$ uniformly on
compact sets, then $\|\phi_0\|_\infty\geq c.$ By assumption (f3)
$\phi_0=\sum_{n=1}^N a_n \frac{\partial w}{\partial x_n}$. On the
other hand for $m=1,\ldots, N$, using \eqref{esti22},  $ 0=\intr
\phi_k(\e_k x+P_1^k)Z_{P_1^k,m}(\e_k x+P_1^k)\to \sum_{n=1}^N a_n
\intr \frac{\partial w}{\partial x_n}(1-\Delta) \frac{\partial
w}{\partial x_m}=a_m\|\frac{\partial w}{\partial
x_1}\|_{H^1(\R^N)}^2$, which implies $a_m=0$, that is $\phi_0=0$.
The contradiction follows.

Hence we have proved  (\ref{cla}), by which  we immediately obtain
\[ \big\|f' (\chi w_{{\bf P }_k}) \phi_k \big\|_{*,{\bf P}_k} =o(1)
\]
and, by (\ref{del1}),
\begin{equation*}
\label{del3} \| \e_k^2\Delta \phi_k -V(x)\phi_k\|_{*,{\bf
P}_k}=o(1)
\end{equation*}
Observe that by \eqref{wdecay}, if we set $\Phi_k(x)=\frac{1}{2}\big(\sum_{i=1}^\ell w_{\epii^k}\big)^{\mu}$, it
follows that,  provided that $\mu$ is chosen sufficiently small, for every $k$:
$$\e_k^2\Delta \Phi_k-V(x)\Phi_k \leq -\frac{\inf_{\R^N}V}{2}\Phi_k
 \hbox{ in }\R^N.$$
 Then one has
\begin{equation*}\e_k^2\Delta (\Phi_k\pm\phi_k)-V(x)(\Phi_k\pm\phi_k)\leq 0\hbox{ in }\R^N.
\end{equation*}
By the comparison principle it follows that $\Phi_k\pm\phi_k\geq
0$.
 Then we have $|\phi_k|\leq\frac{1}{2} \big(\sum_{i=1}^\ell w_{\epii^k}\big)^\mu$,
 by which
$\|\phi_k\|_{*, {\bf P}_k}\leq \frac{1}{2}$, in contradiction with
$\|\phi_k\|_{*, {\bf P}_k}=1$. \hfill$\Box$

\bigskip
Now we are in position to provide the existence of a solution for
the system (\ref{linear}).

\begin{lemma}\label{ex}
For $\e>0$ sufficiently small, for every ${\bf P }\in \Gamma_\e$
and  $ \theta \in L^2(\R^N)$, there exists a unique pair $ (\phi,
\alpha_{in}) $ solving (\ref{linear}). Furthermore
\begin{equation*}
\label{phih} \|\phi\|_{*,{\bf P}} \leq  C\|\theta\|_{*,\bf
P},\quad |\alpha_{in}|\leq C(\e^{1+\sigma}\|\theta\|_{*, {\bf
P}}+\e\|\theta\|_{*, {\bf P}}) .
\end{equation*}
\end{lemma}

\noindent {\textit{Proof.}} The existence follows from Fredholm's
alternative. For every ${\bf P }\in\Gamma_\e$ let us consider
${\cal H}_{\bf P }$ the closed subset of $H^1_V(\R^N)$ defined by
   \[{\cal H}_{\bf P }= \Big\{ \phi \in H^1_V(\R^N) \, \Big|\,
 \Big(\phi,\frac{\partial (\chi w_{P_i})}{\partial x_n} \Big)_\e =0
  \; \;\forall i=1,\ldots, \ell,\;\forall n=1,\ldots N\Big\}.\]
Notice  that, by (\ref{ort}),  $ \phi\in {\cal H}_{\bf P }$ solves
the equation ${\cal L}_{{\bf P
}}[\phi]=\theta+\sum_{i,n}\alpha_{in}Z_{P_i, n}$ if and only  if
   \begin{equation}\label{exit} (\phi, \psi)_\e  -  \intr f'(\chi w_{\epi})\phi\psi dx =
   -\intr \theta
\psi \,dx \;\;  \forall \psi \in {\cal H}_{\bf P }.\end{equation}
Indeed, once we know $\phi$, we can determine the unique
$\alpha_{in}$ from the linear system of equations
\begin{align*}   \intr f'(\chi w_{\epi})\phi\frac{\partial (\chi w_{\epij})}{\partial x_m} \,dx =\intr \theta
   \frac{\partial (\chi w_{\epij})}{\partial x_m}\,dx  +\sum_{i,n}\alpha_{in}
\intr Z_{P_i,n} \frac{\partial (\chi w_{\epij})}{\partial
   x_m}\, dx,\end{align*} for $j=1,\ldots, \ell,\,m=1,\ldots, N, $ which is uniquely solvable according to \eqref{dom}.
   By standard elliptic regularity, $\phi\in
   H^2(\R^N)$.

   Thus it remains to solve (\ref{exit}). According to Riesz's representation theorem, take
   ${\cal K}_{{\bf P }}(\phi)$, $\overline{\theta}\in {\cal H}_{\bf P }$
such that $$({\cal K}_{\bf P }(\phi), \psi)_\e  =- \intr f'(\chi
w_\epi)\phi\psi\, dx \;\;\;\;(\overline{\theta}, \psi)_\e=- \intr
\theta \psi\,dx \;\;\;\;\forall \psi \in {\cal H}_{\bf P }.$$ Then
problem (\ref{exit}) consists in finding $\phi\in {\cal H}_{\bf P
}$ such that
   \begin{equation}
   \label{phhh}
   \phi + {\cal K}_{\bf P } (\phi) = \overline{\theta}.
   \end{equation}
   It is easy to prove that ${\cal K}_{\bf P }$ is a linear compact operator from ${\cal H}_{\bf P }$
    to ${\cal H}_{\bf P }$.
   Using Fredholm's alternatives,  (\ref{phhh}) has a unique  solution for each $\overline{\theta}$, if and only if (\ref{phhh})  has a unique solution
for $\overline{\theta}=0$. Let   $\phi\in {\cal H}_{\bf P }$ be a
solution of $\phi+{\cal K}_{\bf P }(\phi)=0$; then  $\phi$ solves
the system (\ref{linear}) with $\theta=0$ for some
$\alpha_{in}\in\R$. Lemma \ref{apriori} implies $\phi\equiv 0.$
The remaining part of the Lemma follow by Lemma \ref{apriorii} and
Lemma \ref{apriori}.\hfill $\Box$

\subsection{Lyapunov-Schmidt Reduction}
To complete the Lyapunov-Schmidt Reduction, it remains to  prove  Lemma \ref{reg} and Proposition \ref{energy111}.

\noindent{\bf{Proof of Lemma \ref{reg}.}} We write the equation in
(\ref{nonl}) in the following form:
\begin{equation}
\label{33} {\cal L}_{{\bf P }} [\phi] = -{\cal S}_\e [\chi
w_{\epi}] - N_{{\bf P }} [\phi] +\sum_{i,n} \alpha_{in} Z_{P_i, n}
\end{equation}
and use contraction mapping theorem. Here
\begin{equation*}
\label{Mep} {\cal N}_{{\bf P }} [\phi]=  f(\chi w_{\epi} +\phi)-
f(\chi w_{\epi}) - f' (\chi w_{\epi}) \phi   .
\end{equation*}
Consider the metric space  $ {\mathcal B}_{\bf P }=  \{ \phi\in
L^2(\R^N)\, |\,\| \phi \|_{*,{\bf P}}\leq \e ^{\eta} \}$ endowed
with the norm $\|\cdot\|_{*,{\bf P}}$. Given $\phi_1$, $ \phi_2\in
{\mathcal B}_{\bf P }$, by assumption (f1) we have
\begin{equation}\begin{aligned}\label{szz}\|{\cal N}_{{\bf P }}[\phi_1]-{\cal N}_{{\bf P }}[\phi_2]\|_{\ast,{\bf P}}&
\leq C\e ^{\sigma\eta}\|\phi_1-\phi_2\|_{*,{\bf P}}.
\end{aligned}\end{equation}
\noindent For every $\phi\in {\mathcal B}_{\bf P }$  we define
${\mathcal A}_{\bf P } [\phi] \in H^2(\R^N)\cap H^1_V(\R^N) $ to
be the unique solution to the system (\ref{linear}) given by Lemma
\ref{ex} with $ \theta=\theta_{\bf P }[\phi]:= -{\cal S}_{\e}[\chi
w_{{\bf P }}]-{\cal N}_{{\bf P }}[\phi] $.
 By \eqref{szz}, Lemma \ref{error}, Lemma \ref{ex}
$$\| {\mathcal A}_{\bf P }  [\phi]\|_{*,{\bf P}}\leq C \|\theta_{\bf P }[\phi]\|_{*,\bf P}
\leq C(\e^{\beta(\beta+\sigma)}+
    \e^{(1+\sigma)\eta})< \e^{\eta}$$ at least for small $\e$, and hence $ {\mathcal A}_{\bf P }
  [\phi] \in
{\mathcal B}_{\bf P }$.
 Moreover, since
${\mathcal A}_{\bf P }  [\phi_1] -{\mathcal A}_{\bf P }  [\phi_2]$
    solves the system (\ref{linear}) with $\theta=-{\cal N}_{{\bf P }}  [\phi_1] +{\cal N}_{{\bf P }}  [\phi_2]$, by (\ref{szz})  and Lemma \ref{ex}
    we also have that
   \begin{equation*}
   \label{tep2}\begin{aligned}
   \| {\mathcal A}_{\bf P }  [\phi_1] -{\mathcal A}_{\bf P }  [\phi_2]\|_{*,{\bf P}}
 &\leq
   C \| {\cal N}_{{\bf P }}  [\phi_1] -{\cal N}_{{\bf P }}  [\phi_2] \|_{*,{\bf P}}
 <\| \phi_1 -\phi_2 \|_{*,{\bf P}}\;\; \forall \phi_1,\,\phi_2\in
{\mathcal B}_{\bf P },\;\;
   \;\;\forall {\bf P }\in
\Gamma_\e,\end{aligned}
   \end{equation*}
  i.e. the map ${\mathcal A}_{\bf P } $ is a contraction map from  ${\mathcal B}_{\bf P }$ to ${\mathcal B}_{\bf P }$.
   By the contraction mapping theorem, (\ref{nonl}) has a unique solution  $ (\phi_{\bf P }, \alpha_{in}({\bf P } ) )
    \in {\mathcal B}_{\bf P }\times \R^{N\ell}$.

    Finally, by multiplying the equation in (\ref{33}) by $\phi_{\bf P}$
and integrating over $\R^N$ we immediately obtain $(\phi_{\bf
P},\phi_{\bf P})_\e\leq C \e^{N+2\eta}$.
     By Lemma \ref{apriorii}  we get
    $$ \begin{aligned}|\alpha_{in}({\bf P})|&\leq C(\e^{1+\sigma}\|\phi_{\bf P}\|_{*, {\bf P}}+\e\|\theta_{\bf P}[\phi_{\bf P}]\|_{*, {\bf P}})
    \leq C\e^{1+\eta}.\end{aligned}$$

The fact that the map ${\bf P }\in \Gamma_\e \to \phi_{{\bf P
}}\in H^1_V(\R^N)$ is ${\mathcal C}^1$ follows from the Implicit
Function Theorem. See \cite{amma}, for instance.

\bigskip

\noindent {\bf Proof of Proposition \ref{energy111}.} We compute
\begin{align*}I_\e [
\chi w_{\epi} + \phi_{{\bf P }}] &=\frac{1}{2} \intr\big(\e^2
|\nabla (\chi w_{{\bf P }} + \phi_{{\bf P }})|^2 + V(x)(\chi
w_{{\bf P }}
+ \phi_{{\bf P }})^2 \big) dx -\intr F(\chi w_{{\bf P }} + \phi_{{\bf P }}) dx\\
&= I_\e[\chi w_{\bf P }]-\intr {\cal S}_\e [\chi w_{\bf P }]
\phi_{\bf P } dx+ \fr (\phi_{\bf P}, \phi_{\bf P})_\e\\
&\;\;\;\;-\intr \big(F(\chi w_{{\bf P }} + \phi_{{\bf P }})
-F(\chi w_{\bf P })-f(\chi w_{\bf P })\phi_{\bf P }\big) dx.
\end{align*}  By
Lemma \ref{error} we have $|{\cal S}_\e[\chi w_{\bf P }]| \leq
\e^{\eta}\sum_{i=1}^\ell w_\epii^{1-\beta^2}$ for small $\e$,
while $|F(\chi w_{{\bf P }} + \phi_{{\bf P }})-F(\chi w_{\bf P
})-f(\chi w_{\bf P })\phi_{\bf P }|\leq C|\phi_{\bf P }|^{2};$
hence, by using (\ref{phi}) we get
$$I_\e [
\chi w_{{\bf P }} + \phi_{{\bf P }}] =I_\e[\chi w_{\bf P
}]+O(\e^{N+2})$$ uniformly for ${\bf P }\in {\Gamma}_\e$.
\eqref{preespi1} follows from Proposition \ref{energy1}. Next,
denoting by $P_{i,n}$ the $n$-th component of $P_i$,  since
$\frac{\partial w_{\epi}}{\partial P_i^n}=-\tau_{i} \frac{\partial
w_{\epii}}{\partial x_n}$, we  compute
$$\begin{aligned}&\frac{\partial}{\partial P_{i,n}}I_\e [
\chi w_{\epi} + \phi_{{\bf P }}] =-\intr {\cal S}_\e [ \chi
w_{\epi} + \phi_{{\bf P }}]\frac{\partial(\chi w_{\epi}+\phi_{\bf
P})}{\partial P_{i,n}}dx\\ &=\frac{\partial}{\partial P_{i,n}}I_\e
[ \chi w_{\epi}] -\tau_{i}\Big(\phi_{\bf P}, \frac{\partial(\chi
w_{\epii})}{\partial x_n}\Big)_\e-\!\intr  {\cal S}_\e [ \chi
w_{\epi} + \phi_{{\bf P }}]\frac{\partial \phi_{{\bf P
}}}{\partial P_{i,n}}-\!\intr (f(\chi w_\epi+\phi_\epi)-f(\chi w_\epi))\frac{\partial(\chi w_{\epi})}{\partial P_{i,n}}\\
&=\frac{\partial}{\partial P_{i,n}}I_\e [ \chi w_{\epi}]
-\sum_{j,m}\alpha_{jm}({\bf P})\intr Z_{P_j, m} \frac{\partial
\phi_{{\bf P }}}{\partial P_{i,n}}+\tau_{i}\intr (f(\chi
w_\epi+\phi_\epi)-f(\chi w_\epi))\chi\frac{\partial
w_{\epii}}{\partial x_n}.
\end{aligned}
$$
Since $\intr Z_{P_j,m} \phi_\epi\,dx=0$,
 by differentiation
we get\footnote{Observe that  $|\frac{\partial Z_{P_j,m}}{\partial
P_{i,n}}|=\delta_{ij}|(V(x)-\e^2\Delta)\big(\frac{\partial}{\partial
x_n}\big(\chi \frac{\partial w_\epii}{\partial x_m}\big)\big)|\leq
C\e^{-2}w_\epii$ by \eqref{wdecay}.}
\begin{equation}\label{probb}\intr Z_{P_j,m}\frac{\partial \phi_{{\bf P }}}{\partial P_{i,n}} dx=
-\intr\frac{\partial Z_{P_j,m}}{\partial P_{i,n}}   \phi_{{\bf P
}}=O(\e^{N+\eta-2}),\end{equation} by which, using Lemma
\ref{reg},
\begin{equation}\label{biri1}\sum_{j,m}\alpha_{jm}({\bf P})\intr Z_{P_j, m} \frac{\partial \phi_{{\bf
P }}}{\partial P_{i,n}}\, dx=O(\e^{N+2\eta-1}).\end{equation} By
assumption (f1)  we have  $|f(\chi w_{{\bf P }} + \phi_{{\bf P
}})-f(\chi w_{\bf P })-f'(\chi w_{\bf P })\phi_{\bf P }|\leq
C|\phi_{\bf P }|^{1+\sigma};$ consequently
\begin{equation}\label{biri2}\begin{aligned}\intr (f(\chi w_\epi+\phi_\epi)-f(\chi
w_\epi)-f'(\chi w_\epi)\phi_\epi)\chi\frac{\partial
w_{\epii}}{\partial x_n}&= O(\e^{N+\eta(1+\sigma)-1}).
\end{aligned}
\end{equation}
Finally, by \eqref{esti00} and \eqref{esti22},
\begin{equation}\label{biri3}\begin{aligned}\bigg|\intr f'(\chi w_\epi)\phi_\epi\chi\frac{\partial w_{\epii}}{\partial
x_n} \bigg|&=\bigg|\intr \bigg(f'(\chi w_\epi)\chi\frac{\partial
w_{\epii}}{\partial x_n}-Z_{P_i, n}\bigg)\phi_\epi dx\bigg|\\
&\leq C \e^\eta \intr |f'(w_\epi)-f'(w_\epii)|\Big|\frac{\partial
w_{\epii}}{\partial x_n}\Big|dx +C\e^{N+2\beta+\eta-1}\\ &\leq
C\e^\eta\sum_{j\neq i}\intr w_\epij^\sigma \Big|\frac{\partial
w_{\epii}}{\partial x_n}\Big|dx+C\e^{N+2\beta+\eta-1} \leq
C\e^{N+2\beta \sigma+\eta-1}
\end{aligned}\end{equation}
 where in the last inequality we have used \eqref{esti00}.
 Combining \eqref{biri1}-\eqref{biri2}-\eqref{biri3},   we deduce
$$\begin{aligned}\frac{\partial}{\partial P_{i,n}}I_\e [
\chi w_{\epi} + \phi_{{\bf P }}] &=\frac{\partial}{\partial
P_{i,n}}I_\e [ \chi w_{\epi} ] + O(\e^{N+\beta(1+\sigma)^2-1})
\end{aligned}
$$
uniformly for ${\bf P}\in\Gamma_\e$. By applying Proposition
\ref{energy1} we obtain \eqref{preespi2}, using that
 $\beta(1+\sigma)^2-1>\beta$ thanks to assumption  (f1) if $\beta$ is close to 1.
\hfill $\Box$

\bigskip
\medskip

\small{\noindent{\bf Acknowledgments.}  This paper was begun  while the first author was visiting the University
of Granada in February 2009. She gratefully acknowledges the  Department of Mathematical Analysis for their kind
hospitality.}

\end{document}